\documentclass[12pt,a4paper,reqno]{amsart}

\usepackage{amsfonts}
\usepackage{amsmath}
\usepackage{amssymb}
\usepackage{amsthm}
\usepackage{amsxtra}
\usepackage{bbm}
\usepackage{braket}
\usepackage{comment}
\usepackage{extarrows}

\usepackage{graphicx}
\usepackage{latexsym}

\usepackage{mathrsfs}
\usepackage{yhmath}

\usepackage{nicematrix}
\usepackage{mathtools}

\usepackage{float}
\usepackage{booktabs}
\usepackage{verbatim}
\usepackage{xcolor}

\usepackage[pdfborder={0 0 0}]{hyperref} 
\hypersetup{breaklinks, raiselinks}

\usepackage{bm}

\usepackage{cleveref}

\usepackage{a4wide}
\usepackage{fullpage}

\theoremstyle{plain}
\newtheorem{Theorem}{Theorem}[section]

\newtheorem{Lemma}[Theorem]{Lemma}
\newtheorem{Corollary}[Theorem]{Corollary}
\newtheorem{Proposition}[Theorem]{Proposition}

\theoremstyle{definition}

\newtheorem{Assumptions and Discussion}[Theorem]{Assumptions and Discussion}

\newtheorem{Definition}[Theorem]{Definition}

\theoremstyle{remark}

\newtheorem*{acknowledgment*}{Acknowledgment}

\numberwithin{equation}{section}

\usepackage[shortlabels]{enumitem}
\SetEnumerateShortLabel{a}{\textup{(\alph*)}}
\SetEnumerateShortLabel{A}{\textup{(\Alph*)}}
\SetEnumerateShortLabel{1}{\textup{(\arabic*)}}
\SetEnumerateShortLabel{i}{\textup{(\roman*)}}
\SetEnumerateShortLabel{I}{\textup{(\Roman*)}}

\begin{document}

\title{Curvature estimates for stable hypersurfaces of the Euler-Dierkes-Huisken functional}

\author[Hongbin Cui, Xiaowei Xu, Shoulan Yi]{Hongbin Cui \quad Xiaowei Xu\quad  Shoulan Yi}

\address{Hongbin Cui, School of Mathematical Sciences, University of Science and Technology of China, Hefei, Anhui, 230026, P.R.~China}
\email{cuihongbin@ustc.edu.cn}

\address{Xiaowei Xu (The corresponding author), School of Mathematical Sciences, University of Science and Technology of China, Hefei, Anhui, 230026, P.R.~China}
\email{xwxu09@ustc.edu.cn}

\address{Shoulan Yi, School of Mathematical Sciences, University of Science and Technology of China, Hefei, Anhui, 230026, P.R.~China}
\email{yishoulan@mail.ustc.edu.cn}

\begin{abstract}
A hypersurface $\Sigma \subset\mathbb{R}^{n+1}$ is called EDH-stable if it is a critical point of the Euler-Dierkes-Huisken functional $\mathcal{E}_\alpha(\Sigma)=\int_\Sigma|x|^\alpha\, d\mathcal{H}^n$ with nonnegative second variation, which was introduced by Dierkes and Huisken in \cite{[DH-24]}. In this paper, we first establish a Bellettini-type $\epsilon$-regularity theorem and use it, together with tangent cone analysis, to derive curvature estimates for EDH-stable hypersurfaces with $2\leq n \leq 6$.

\end{abstract}
\maketitle
\tableofcontents

\section{Introduction}
In 1744, L. Euler considered the least moment of inertia for planar curves with fixed endpoints and introduced the functional
\begin{equation*}
\mathcal{E}_2(x(t))=\int_a^b|x(t)|^2\, ds,
\end{equation*}
where $x(t)=(x_1(t),x_2(t))\in\mathbb{R}^2$, $t\in[a,b]$, and $ds=|x'(t)|dt$. M. Mason (\cite{[M-1906]}), L. Tonelli (\cite{[T-1921]}), C. Caratheodory (\cite{[C-1935]}) studied regular extremals of a more general functional
\begin{equation*}
\mathcal{E}_\alpha(x(t))=\int_a^b|x(t)|^\alpha\, ds,
\end{equation*}
where $\alpha>0$. Recently, U. Dierkes and G. Huisken (\cite{[DH-24]}) considered the $n$-dimensional analogue functional
\begin{equation}\label{A-1}
\mathcal{E}_\alpha(\Sigma):=\int_\Sigma|x|^\alpha\, d\mathcal{H}^n
\end{equation}
for hypersurfaces in $\mathbb{R}^{n+1}$, where $x$ is the position of a point in $\mathbb{R}^{n+1}$, $\alpha\in \mathbb{R}$, and $\mathcal{H}^n$ is the $n$-dimensional Hausdorff measure on $\mathbb{R}^{n+1}$. We call (\ref{A-1}) the \emph{Euler-Dierkes-Huisken functional} (EDH-functional) in the sequel. A hypersurface $\Sigma \subset\mathbb{R}^{n+1}$ is called \emph{EDH-stable} if it is a critical point of the EDH-functional with nonnegative second variation. When $\alpha=0$, the EDH-functional is simply​ the area functional, and a critical point of the latter​ is a minimal surface.

The resolutions of the graphical Bernstein problem and the stable
Bernstein problem have played an important role in the development
of the theory of minimal surfaces.
In 1975, R. Schoen, L. Simon and S.T. Yau (\cite{[SSY-1975]}) proved that any properly immersed two-sided stable minimal hypersurface $\Sigma$ in $\mathbb{R}^{n+1}$ with $n\leq 5$ and Euclidean volume growth, is a union of affine hyperplanes. N. Wickramasekera \cite{[W-2008]} proved the case $n=6$ under a multiplicity-two
condition and an explicit volume-growth condition. R. Schoen and L.Simon (\cite{[SS-1981]}) proved that the results hold for $n\leq 6$ if $\Sigma$ is properly embedded. Recently, C. Bellettini \cite{[B-2025]} extended the Schoen--Simon--Yau curvature estimates and the associated Bernstein-type theorems to immersed hypersurfaces
in $\mathbb{R}^{n+1}$ with $n \leq 6$ via an $\epsilon$-regularity theorem and
De Giorgi iteration.

The Bernstein problem for hypersurfaces with respect to weighted area functionals has also drawn considerable attention in geometry and analysis. U. Dierkes (\cite{[D-1993]}, \cite{[D-1995]}) established Bernstein-type theorems for stable graphs and area-minimizing hypersurfaces associated with the weighted area functional, where the weight is a power of the last coordinate. L. Wang (\cite{[W-2011]}), Q. Guang and J. Zhu (\cite{[GZ-2017]}) obtained  Bernstein theorems for graphic self-shrinkers. C. Mooney and Y. Yang (\cite{[MY-2024]}) have completed the solution to the anisotropic Bernstein problem, in which the weight function depends on the normal vector. For results on Bernstein-type theorems for other nonlinear geometric PDEs, we refer the reader to \cite{[M-2024]}.

In \cite{[DH-24]}, U. Dierkes and G. Huisken established fundamental properties of EDH-stationary hypersurfaces and investigated the stability and area-minimizing properties​ of hypercones associated with the functional (\ref{A-1}). H.B. Cui and X.W. Xu (\cite{[CX-2025]}) extended the theory of Dierkes and Huisken for​ higher codimension surfaces in Euclidean space. U. Dierkes and  R. López (\cite{[DL-2026]}) proved that an axially symmetric EDH-stationary surface in $\mathbb{R}^3\setminus \{0\}$  with $\alpha>0$ is an entire graph defined on the plane perpendicular to the axis of symmetry. Indeed, one can check that the example constructed by U. Dierkes and  R. López (\cite{[DL-2026]}) is EDH-stable. For the classification of EDH-stationary surfaces in $\mathbb{R}^3$ under various geometric assumptions, we refer the reader to \cite{[DL-2026]}, \cite{[L-2025-1]}, and \cite{[L-2025-2]}.

Regarding the Bernstein-type theorem for $\alpha$-stationary hypersurfaces, H.B. Cui, J.H. Li and X.W. Xu (\cite{[CLX-2026]}) proved that any smooth, complete, connected, embedded EDH-stationary hypersurface $\Sigma\subset \mathbb{R}^{n+1}$ with $\alpha>0$ passing through the origin is a linear hyperplane. The proof relies primarily on the analysis of equation (\ref{A-2}), which differs from classical curvature estimates. However, the examples of Dierkes and López \cite{[DL-2026]} show that there exist nonplanar EDH-stable hypersurfaces in $\mathbb{R}^3$ that do not pass through the origin. Inspired by Bellettini's work \cite{[B-2025]}, in this paper we establish curvature estimates for EDH-stable hypersurfaces in $\mathbb{R}^{n+1}$ with $n \leq 6$.

 We introduce the function $\rho:=\alpha\log|x|$ on $\mathbb{R}^{n+1}\setminus\{0\}$, and denote its differential and Hessian by $D\rho$, $D^2\rho$, respectively. Then, the Euler-Lagrange equation of the EDH-functional is
\begin{equation}\label{A-2}
H=\rho_\nu=\alpha |x|^{-2}\langle x,\nu\rangle,
\end{equation}
and the stability is equivalent to
\begin{equation}\label{A-3}
\int_{\Sigma}(|A|^2-\rho_{\nu\nu})\varphi^2 |x|^\alpha d\mathcal{H}^n\leq\int_{\Sigma}|\nabla\varphi|^2 |x|^\alpha d\mathcal{H}^n,
\end{equation}
for all $\varphi\in C_c^\infty(\Sigma\setminus\{0\})$, where $H$ is the mean curvature of $\Sigma$, and $\nu$ is the unit normal vector field of $\Sigma$, and $\rho_\nu=D\rho(\nu)$, $\rho_{\nu\nu}=D^2\rho(\nu,\nu)$.

Let $p\in\Sigma$, and let $B_{2R}(p)\subset \mathbb{R}^{n+1}\setminus\{0\}$ be relatively compact, and define
\begin{equation}\label{A-4}
a:=|p|-2R>0,\quad
\sigma:=a^{-1},\quad
\lambda:=R\sigma=\frac{R}{|p|-2R},
\end{equation}
and
\begin{equation}\label{A-5}
\mathfrak{C}(p,2R):=R^{2-n}\int_{\Sigma\cap B_{2R}(p)}(|A|+\sigma)^2\,d\mu,
\end{equation}
where $d\mu:=(\frac{|x|}{|p|})^\alpha d\mathcal{H}^n$.

\vspace{0.2cm}

We first have the following $\epsilon$-regularity theorem.

\begin{Theorem}\label{Theorem1.1}
Let $2\leq n \leq 6$ and $\alpha \in \mathbb{R}$, and let $\Sigma\rightarrow\mathbb{R}^{n+1}$ be a 
smooth, two-sided, EDH-stable hypersurface properly immersed in $\mathbb{R}^{n+1}$ with $p\in\Sigma$ and $\partial\Sigma\cap B_{2R}(p)=\emptyset$. There exist constants $\epsilon(n,\alpha)>0$ and $C(n,\alpha)>0$ depending only on $n$ and $\alpha$, such that if
\begin{equation}\label{A-6}
\lambda\leq 1,\quad \mathfrak{C}(p,2R)+\lambda^2\leq\epsilon(n,\alpha),
\end{equation}
then
\begin{equation}\label{A-7}
R\sup\limits_{\Sigma\cap B_{\frac{R}{2}}(p)}|A|\leq C(n,\alpha)\,(\mathfrak{C}(p,2R)^{\frac{1}{2}}+\lambda).
\end{equation}
\end{Theorem}

Theorem \ref{Theorem1.1} admits a more general extension, which will be used frequently in cone approximation.
The proof is similar to that of Theorem \ref{Theorem1.1}.

\begin{Corollary}\label{Corollary1.2}
Let $\Sigma$ be a two-sided, properly immersed stable hypersurface for a more general weighted functional $\int_\Sigma e^\rho d\mathcal{H}^n$ \emph{(}or $\rho$-weighted stable for short\emph{)} in $B_{2R}(p)$ with $\partial\Sigma\cap B_{2R}(p)=\emptyset$, and set $\tilde{\rho}=\rho-\rho(p)$. Assume that 
\begin{equation}\label{A-8}
\sup\limits_{B_{2R}(p)}R^m|D^m\rho|\leq c_0\lambda^m, \quad 1\leq m\leq 3,
\quad
\sup\limits_{B_{2R}(p)}|\tilde{\rho}|\leq c_0\lambda,
\end{equation}
where $0\leq \lambda \leq 1$ and $c_0\geq 1$. Set 
$$\mathfrak{C}_\rho(p,2R):=R^{2-n}\int_{\Sigma\cap B_{2R}(p)}(|A|+\frac{\lambda}{R})^2\,d\mu_{\tilde{\rho}},\quad d\mu_{\tilde{\rho}}=e^{\tilde{\rho}}\,d\||\Sigma|\|.$$
Then, there are $\epsilon_0(n,c_0)>0$ and $C(n,c_0)<\infty$ such that if $\mathfrak{C}_\rho(p,2R)+\lambda^2\leq \epsilon_0(n,c_0)$, we have
\begin{equation}\label{A-9}
R\sup\limits_{\Sigma\cap B_{\frac{R}{2}}(p)}|A|\leq C(n,c_0)(\mathfrak{C}_\rho(p,2R)^{\frac{1}{2}}+\lambda).
\end{equation}
\end{Corollary}

Using Allard's compactness theorem, tangent cone analysis, and Federer's
dimension reduction argument, we establish the following local curvature estimates.

\begin{Theorem}\label{Theorem1.3}
Fix $2\leq n\leq 6$ and $\alpha\in\mathbb{R}$. Let $p_0\in\Sigma$, $R>0$, and set $r_0:=|p_0|$. Assume that

\emph{(1)} $8R\leq r_0$;

\emph{(2)} $\Sigma\cap B_{4R}(p_0)$ is smooth, two-sided, EDH-stable, without boundary in $B_{4R}(p_0)$, and the immersion is proper over $\overline{B}_{4R}(p_0)$;

\emph{(3)} for some $\Lambda<\infty$,
\begin{equation}\label{A-10}
\frac{1}{r^\alpha_0(4R)^n}\int_{\Sigma\cap B_{4R}(p_0)}|x|^\alpha\, d\mathcal{H}^n\leq \Lambda.
\end{equation}
Then there exists a constant $C(n,\alpha,\Lambda)>0$ such that
\begin{equation}\label{A-11}
\sup\limits_{\Sigma\cap B_R(p_0)}|A_\Sigma|\leq \frac{C(n,\alpha,\Lambda)}{R}.
\end{equation}
\end{Theorem}

\emph{Remark}. An alternative proof of Theorem \ref{Theorem1.3} follows directly from a
Bernstein-type theorem (Theorem 1 in \cite{[B-2025]}) and blow-up analysis.
Here we instead develop Bellettini's theory further and extend it to
the setting of EDH-stable hypersurfaces. The estimates in Theorem \ref{Theorem1.3} are local, which hold away from the origin. Since $|D^m\rho|=O(r^{-m})$ as $r\rightarrow 0$, the metric fails to have bounded geometry near the origin. Consequently, the estimates in Lemma \ref{Lemma4.1} and Lemma \ref{Lemma4.4} are too weak to support a small-energy in cone approximation argument. 

The paper is organized as follows. In Section 2, we derive the weighted Simons
inequality and a refined Kato inequality. In Section 3, we establish a
Caccioppoli inequality and prove the $\epsilon$-regularity of the second
fundamental form of EDH-stable hypersurfaces. In Section 4, we carry out the
geometric blow-up argument for EDH-stable hypersurfaces, i.e., the cone
approximation of admissible sequences. In Section 5, we perform the tangent
cone analysis of the integral varifold arising from an admissible sequence. In Section 6, we
give the proof of our main theorem.
\vspace{0.5cm}

\section{Weighted Simons inequality and Kato inequality}
The goal of this section is to derive the weighted Simons inequality and a refined Kato inequality, which arise as a special case of the Simons-type formula for $f$-minimal hypersurfaces due to Cheng, Mejia, and Zhou (\cite{[CMZ-2015]}) and as a natural application of the trace-free fully symmetric three-tensor, respectively.

Let $x:\Sigma^n\rightarrow \mathbb{R}^{n+1}$ be a smooth immersion. Choose a local orthonormal frame $e_1,\ldots,e_n,\nu$ such that $e_1,\ldots, e_n$ are tangent to $\Sigma$ and $\nu$ is normal to $\Sigma$. The coefficients of the second fundamental form $A$ with respect to the chosen frame are $h_{ij}=\langle D_{e_i}e_j,\nu\rangle$. By the Codazzi equation and Gauss equation, we obtain
\begin{equation}\label{B-1}
\Delta h_{ij}=H_{ij}+Hh_{ik}h_{kj}-|A|^2h_{ij},
\end{equation}
where $H$ is the mean curvature. On the other hand, using the Euler-Lagrange equation $H=\rho_\nu$, we have 
\begin{equation}\label{B-2}
H_i=\rho_{i\nu}-h_{ik}\rho_k,
\end{equation}
\begin{equation}\label{B-3}
H_{ij}=\rho_{\nu ij}+h_{ij}\rho_{\nu\nu}-h_{ijk}\rho_k-h_{ik}\rho_{kj}-\rho_{ik}h_{kj}-Hh_{ik}h_{kj},
\end{equation}
where the Codazzi equation has been used. We adopt the notation that $\rho_i=D\rho(e_i)$, $\rho_{ij}=D^2\rho(e_i,e_j)$, $\rho_{\nu ij}=D^3\rho(\nu,e_i,e_j)$, and so on. Substituting (\ref{B-3}) into (\ref{B-1}) and using the drift Laplacian $\Delta_\rho =\Delta +\langle\nabla\rho, \nabla\rangle$, we have
\begin{equation}\label{B-4}
\Delta_\rho h_{ij}=-|A|^2h_{ij}+\rho_{\nu ij}+h_{ij}\rho_{\nu\nu}-h_{ik}\rho_{kj}-\rho_{ik}h_{kj}.
\end{equation}
Then
\begin{equation*}
\frac{1}{2}\Delta_\rho |A|^2=|\nabla A|^2+h_{ij}\Delta_\rho h_{ij}
\end{equation*}
gives the weighted Simons identity
\begin{equation}\label{B-5}
\frac{1}{2}\Delta_\rho |A|^2=|\nabla A|^2-|A|^4+|A|^2\rho_{\nu\nu}+h_{ij}\rho_{\nu ij}-2h_{ij}h_{ik}\rho_{kj}.
\end{equation}

The following refined Kato inequality is a standard consequence of the orthogonal trace-free decomposition of the Codazzi tensor. We include the proof only to keep track of the dimensional constants.

\begin{Proposition}\label{Proposition2.1}
At the point where $|A|>0$, the following pointwise estimate holds
\begin{equation}\label{B-6}
|\nabla A|^2-|\nabla|A||^2\geq \frac{1}{n} |\nabla|A||^2-C_n|\nabla H|^2,
\end{equation}
where $C_n=\frac{(n+1)(2+\sqrt{n})^2}{n(n+2)}$.
\end{Proposition}

\emph{Proof}. Set $S(\nabla H)_{ijk}:=\frac{1}{n+2}(H_i\delta_{jk}+H_j\delta_{ik}+H_k\delta_{ij})$ and  $\overset{\circ}{h}_{ijk}:=h_{ijk}-S(\nabla H)_{ijk}$. Then $\nabla A=\overset{\circ}{\nabla A}+S(\nabla H)$ is an orthogonal decomposition, and hence
\begin{equation}\label{B-7}
|\nabla A|^2=|\overset{\circ}{\nabla A}|^2+|S(\nabla H)|^2=|\overset{\circ}{\nabla A}|^2+\frac{3}{n+2}|\nabla H|^2,
\end{equation} 
where $\mathring{\nabla A}$ is the trace-free part of $\nabla A$. 
Noting that $|A||A|_k = h_{ij}h_{ijk}$ and using
$\mathring{h}_{ij}=h_{ij}-\frac{H}{n}\delta_{ij}$ and the fact that $\sum\limits_{i}\overset{\circ}{h}_{iik}=0$, we compute
\begin{eqnarray}\label{B-8}
|\nabla |A||&\leq& |A|^{-1}|h_{ij}\overset{\circ}{h}_{ijk}|+|A|^{-1}|h_{ij}S(\nabla H)_{ijk}|\nonumber\\
&=&|A|^{-1}|\overset{\circ}{h}_{ij}\overset{\circ}{h}_{ijk}|+|A|^{-1}|h_{ij}S(\nabla H)_{ijk}|.
\end{eqnarray}

We now estimate the last two terms on the right-hand side of \eqref{B-8}. Given a unit covector $\zeta$, set
$P_{ijk}:=\frac{1}{3}(\mathring{h}_{ij}\zeta_k+\mathring{h}_{ik}\zeta_j+\mathring{h}_{jk}\zeta_i)$. Let $\mathring{P}$ denote the trace-free part of $P$. Then, we obtain
\begin{equation}\label{B-9}
|\mathring{P}|^2=\frac{1}{3}|\mathring{h}_{ij}|^2|\zeta|^2+\frac{2n}{3(n+2)}|\mathring{h}_{jk}\zeta_j|^2.
\end{equation}
Since $\mathring{h}_{ij}$ is trace-free, we have $|\mathring{h}_{jk}\zeta_j|^2\leq \frac{n-1}{n} |\mathring{h}_{jk}|^2|\zeta|^2$. Hence, using \eqref{B-9}, one gets $|\mathring{P}|^2\leq \frac{n}{n+2}|\mathring{h}_{ij}|^2|\zeta|^2$. Observing that $\langle\mathring{\nabla A},P\rangle=\mathring{h}_{ij}\mathring{h}_{ijk}\zeta_k$, we deduce
\begin{equation}\label{B-10}
|\mathring{h}_{ij}\mathring{h}_{ijk}\zeta_k|\leq \sqrt{\frac{n}{n+2}}|\mathring{\nabla A}||\mathring{h}_{ij}||\zeta|.
\end{equation}
Taking the supremum over all unit covector $\zeta$, we obtain from \eqref{B-10} that
\begin{equation}\label{B-11}
|\mathring{h}_{ij}\mathring{h}_{ijk}|\leq\sqrt{\frac{n}{n+2}} |A||\mathring{\nabla A}|.
\end{equation}
On the other hand, a straightforward calculation and $|H|\leq\sqrt{n}|A|$ yield
\begin{equation}\label{B-12}
|h_{ij}S(\nabla H)_{ijk}|\leq \frac{2+\sqrt{n}}{n+2}|A||\nabla H|.
\end{equation}
Substituting \eqref{B-11} and \eqref{B-12} into \eqref{B-8} and applying the elementary inequality $(a+b)^2\leq (1+\epsilon)a^2+(1+\epsilon^{-1})b^2$ with $\epsilon=(1+n)^{-1}$, we obtain
\begin{equation}\label{B-13}
|\mathring{\nabla A}|^2\geq (1+\frac{1}{n})|\nabla|A||^2-C_n|\nabla H|^2,
\end{equation}
where $C_n=\frac{(n+1)(2+\sqrt{n})^2}{n(n+2)}$. Thus the desired assertion follows by combining \eqref{B-7} and \eqref{B-13}.

\hfill$\Box$

\section{Bellettini-type $\epsilon$-regularity theorem}
We first derive the Caccioppoli inequality, and then use it to establish the $\epsilon$-regularity of the second fundamental form​ of EDH-stable hypersurfaces.

\subsection{Caccioppoli inequality}

For notational convenience, we define 
\begin{equation*}
u:=|A|,\quad\quad w:=(u-k)^+
\end{equation*}
and
\begin{equation*}
 E:=|\nabla A|^2-|\nabla u|^2, \quad\quad F:=\rho_{\nu\nu}u^2+\rho_{\nu ij}h_{ij}-2h_{ij}h_{jk}\rho_{ki}.
\end{equation*}
\begin{Proposition}\label{Proposition3.1}
Let $k\geq 0$ and $\xi\in C^1_c(B_{2R}(p))$. 
Then we have
\begin{eqnarray}\label{C-1}
\int_{\{u>k\}}(1-\frac{k}{u})E\xi^2d\mu&\leq& \int_{\{u>k\}} w^2|\nabla\xi|^2 d\mu+k\int_{\{u>k\}}u^2w\xi^2 d\mu\nonumber\\
&{}&-\int_{\{u>k\}}(1-\frac{k}{u})F\xi^2 d\mu+\int_{\{u>k\}}\rho_{\nu\nu}w^2\xi^2 d\mu.
\end{eqnarray}
\end{Proposition}
	​
\emph{Proof}. Since $p\in\Sigma $ is fixed, the stability inequality implies
\begin{equation}\label{C-2}
\int_\Sigma (u^2-\rho_{\nu\nu})\varphi^2\, d\mu\leq\int_{\Sigma}|\nabla\varphi|^2 d\mu,
\end{equation}
for all $\varphi\in C_c^\infty (\Sigma\setminus\{0\})$. Taking $\varphi=w\xi$ in \eqref{C-2}, we have
\begin{equation}\label{C-3}
\int_\Sigma u^2w^2\xi^2 d\mu\leq\int_\Sigma|\nabla(w\xi)|^2 d\mu+\int_\Sigma\rho_{\nu\nu}w^2\xi^2d\mu.
\end{equation}
Using integration by parts, the first term on right-hand side of \eqref{C-3} becomes 
\begin{equation}\label{C-4}
\int_\Sigma|\nabla(w\xi)|^2 d\mu=\int_{\{u>k\}}|\nabla u|^2\xi^2 d\mu-\frac{1}{2}\int_{\{u>k\}}\Delta_\rho w^2\xi^2 d\mu+\int_{\{u>k\}}w^2|\nabla\xi|^2d\mu.
\end{equation}
Using $w=u-k$ on $\{u>k\}$, together with the Simons identity 
$$\frac{1}{2}\Delta_\rho u^2-|\nabla u|^2=u\Delta_\rho u=E-u^4+F,$$
we obtain
\begin{equation}\label{C-5}
\frac{1}{2}\int_{\{u>k\}}\Delta_\rho w^2\xi^2 d\mu=\int_{\{u>k\}}(E-u^4+F+|\nabla u|^2-k\Delta_\rho u)\xi^2\, d\mu.
\end{equation}
Substituting \eqref{C-5} into \eqref{C-4} and using the Simons identity again, we obtain from \eqref{C-3} 
\begin{eqnarray}\label{C-6}
\int_{\{u>k\}} u^2w^2\xi^2 d\mu&\leq&\int_{\{u>k\}}w^2|\nabla\xi|^2d\mu-\int_{\{u>k\}}(1-\frac{k}{u})(E+F)\xi^2d\mu\nonumber\\
&{}&+\int_{\{u>k\}}u^4\xi^2d\mu-k\int_{\{u>k\}}u^3\xi^2d\mu+\int_{\{u>k\}}\rho_{\nu\nu}w^2\xi^2d\mu.
\end{eqnarray}
Thus, the Caccioppoli inequality \eqref{C-1} follows by substituting $u^2w^2=u^4-2ku^3+k^2u^2$ on the set $\{u>k\}$ into \eqref{C-6}.

\hfill$\Box$

Observe that the radial density $\rho$ satisfies
\begin{equation}\label{C-7}
D\rho=\alpha\frac{x}{|x|^2}, \quad D^2\rho=\alpha(\frac{I}{|x|^2}-2\,\frac{x\otimes x}{|x|^4}),
\quad |D^m\rho|\leq C_{m,\alpha}|x|^{-m}.
\end{equation}
In particular, on any region where $|x|\geq a>0$ and $\sigma=a^{-1}$,
\begin{equation}\label{C-8}
|D\rho|\leq C_\alpha\sigma,\quad
|D^2\rho|\leq  C_\alpha\sigma^2,\quad
|D^3\rho|\leq C_\alpha\sigma^3.
\end{equation}
Consequently, on the ball $B_{2R}(p)$, estimate \eqref{C-8} yields
\begin{equation}\label{C-9}
|\nabla H|^2\leq C_\alpha(\sigma^4+\sigma^2 u^2),\quad
|F|\leq C_{n,\alpha}(\sigma^2u^2+\sigma^3 u).
\end{equation}
Applying the refined Kato inequality \eqref{B-6}, together with \eqref{C-8} and \eqref{C-9}, we obtain
\begin{eqnarray}\label{C-10}
\frac{1}{n}\int_{\{u>k\}}(1-\frac{k}{u})|\nabla u|^2\xi^2 d\mu&\leq&\int_{\{u>k\}}w^2|\nabla \xi|^2d\mu +
k\int_{\{u>k\}}u^2w \xi^2d\mu\nonumber\\
&{}&+C_{n,\alpha}\int_{\{u>k\}}(\sigma^2w^2+\sigma^2u^2+\sigma^3u+\sigma^4)\xi^2d\mu.\hspace{0.8cm}
\end{eqnarray}
Since $u=w+k$ on the set $\{u>k\}$ and $ku^2w\leq (k+\sigma)w^3+2(k+\sigma)^2w^2+(k+\sigma)^3w$, \eqref{C-10}
yeilds
\begin{eqnarray}\label{C-11}
\frac{1}{n}\int_{\{u>k\}}(1-\frac{k}{u})|\nabla u|^2\xi^2 d\mu&\leq&\int_{\{u>k\}}w^2|\nabla \xi|^2d\mu \nonumber\\
&{}&+C_{n,\alpha}\sum\limits_{j=0}^3(k+\sigma)^{4-j}\int_{\{u>k\}}w^j\xi^2d\mu.
\end{eqnarray}
Since the polynomial term $ku^2w$ vanishes at the initial level $k=0$, estimate \eqref{C-10} yields
\begin{equation}\label{C-12}
\int_\Sigma |\nabla u|^2\xi^2 d\mu\leq C_{n,\alpha}\int_\Sigma(u+\sigma)^2(|\nabla \xi|^2+\sigma^2 \xi^2)d\mu.
\end{equation}

\subsection{Proof of Theorem 1.1} We recall the Michael-Simon inequality (see \cite{[MS-1973]}) for hypersurfaces in $\mathbb{R}^{n+1}$: for any nonegative $f\in C^1_c(\Sigma)$,
\begin{equation}\label{C-13}
\big(\int_\Sigma f^{\frac{n}{n-1}}\, d\mathcal{H}^n\big)^{\frac{n-1}{n}}\leq C_n\int_\Sigma(|\nabla f|+|H|f)\, d\mathcal{H}^n, \quad n\geq 2,
\end{equation}
which implies
\begin{equation}\label{C-14}
\big(\int_\Sigma f^{\frac{2n}{n-2}}\,d\mathcal{H}^n\big)^{\frac{n-2}{n}}\leq C_n \int_\Sigma(|\nabla f|^2+H^2f^2)\,d\mathcal{H}^n,\quad n\geq 3.
\end{equation}
To proceed, we require the following weighted version of the Michael–Simon inequality.

\begin{Lemma}\label{Lemma3.2}
With the notations in \emph{\eqref{A-4}}, for any nonnegative $f\in C_c^1(\Sigma \cap B_{2R}(p))$, if $\lambda\leq 1$, we have
\begin{equation}\label{C-15}
\big(\int_\Sigma |f|^4\,d\mu\big)^{1/2}\leq C_{2,\alpha}\big(\int_\Sigma f^2\,d\mu\big)^{1/2}
\big(\int_\Sigma (|\nabla f|^2+\sigma^2f^2)\,d\mu\big)^{1/2},\quad n=2,
\end{equation} 
and
\begin{equation}\label{C-16}
\big(\int_\Sigma f^{\frac{2n}{n-2}}\,d\mu\big)^{\frac{n-2}{n}}\leq C_{n,\alpha} \int_\Sigma(|\nabla f|^2+\sigma^2f^2)\,\,d\mu,\quad n\geq 3.
\end{equation}
\end{Lemma}

\emph{Proof}. Note that $|p|=a(1+2\lambda)$ and $\lambda\leq 1$. On the set $B_{2R}(p)$, we have
\begin{equation}\label{C-17}
3^{-|\alpha|}\leq (\frac{|x|}{|p|})^\alpha \leq 3^{|\alpha|},\quad |H|\leq |\alpha| \sigma.
\end{equation}
For $n=2$, using $f^2$ instead of $f$ in \eqref{C-13} with $f\in C_c^1( B_{2R}(p))$, the Cauchy-Schwarz inequality implies
\begin{eqnarray}\label{C-18}
\big(\int_\Sigma f^4\,d\mathcal{H}^2 \big)^{1/2}&\leq& C_2\big(2\int_\Sigma f|\nabla f|\, d\mathcal{H}^2 +
\int_\Sigma |H|f^2\, d\mathcal{H}^2\big)\nonumber\\
&\leq&\sqrt{5}C_2\big(\int_\Sigma f^2 \,d\mathcal{H}^2\big)^{1/2}\big(\int_\Sigma (|\nabla f|^2+H^2f^2)\,d\mathcal{H}^2\big)^{1/2}.
\end{eqnarray}
It follows from \eqref{C-17} and \eqref{C-18} that
\begin{eqnarray*}
\big(\int_\Sigma f^4\,d\mu\big)^{1/2}&\leq&3^{|\alpha|/2}\big(\int_\Sigma f^4 \,d\mathcal{H}^2\big)^{1/2}\\
&\leq & C_{2,\alpha} \big(\int_\Sigma f^2 \,d\mu\big)^{1/2}\big(\int_\Sigma (|\nabla f|^2+\sigma^2f^2)\,d\mu \big)^{1/2}.
\end{eqnarray*}
This completes the proof of \eqref{C-15}; the proof of \eqref{C-16} is similar.

\hfill$\Box$

An immediate corollary of Lemma \ref{Lemma3.2} is that
\begin{Corollary}\label{Corollary3.3}
For $n\geq 3$, we have
\begin{equation}\label{C-19}
\big(\int_{\Sigma\cap B_R(p)} u^{\frac{2n}{n-2}}\,d\mu\big)^{\frac{n-2}{n}}\leq C'_{n,\alpha} R^{n-4}\mathfrak{C}(p,2R).
\end{equation}
\end{Corollary}

\emph{Proof}. Choose $\zeta\in C^1_c(B_{2R}(p))$ such that $0\leq \zeta\leq 1$, $\zeta=1$ on $B_{R}(p)$ and $|\nabla \zeta|\leq \frac{2}{R}$. Substituting  $f=u\zeta$ in \eqref{C-16} and applying the Cauchy-Schwarz inequality and \eqref{C-12}, we obtain 
\begin{eqnarray*}
\big(\int_{\Sigma\cap B_R(p)} u^{\frac{2n}{n-2}}\,d\mu\big)^{\frac{n-2}{n}}&\leq& C_{n,\alpha}
\int_{\Sigma\cap B_{2R}(p)} (2u^2|\nabla\zeta|^2+2|\nabla u|^2\zeta^2+\sigma^2u^2\zeta^2)\,d\mu\\
&\leq& C'_{n,\alpha}R^{-2}\int_{\Sigma\cap B_{2R}(p)}(u+\sigma)^2d\mu\\
&=&C'_{n,\alpha}R^{n-4}\mathfrak{C}(p,2R),
\end{eqnarray*}
where we have used the fact $\sigma=\frac{\lambda}{R}\leq \frac{1}{R}$ and $u\leq u+\sigma$.

\hfill$\Box$

We now use De Giorgi's iteration method to prove the $\epsilon$-regularity theorem. We only treat the case $3\leq n\leq 6$, using \eqref{C-16}. The proof for the case $n=2$ is similar and follows from \eqref{C-15}. Let $d>0$ be a dimensionless number to be chosen later and define
\begin{equation*}
r_\ell:=R(\frac{1}{2}+2^{-\ell}),\quad k_\ell:=\frac{d}{R}(1-2^{1-\ell}),\quad w_\ell:=(u-k_\ell)^+,
\end{equation*}
and 
\begin{equation*}
\Omega_\ell:=\{u>k_\ell\}\cap B_{r_\ell}(p),\quad S_\ell:=\int_{B_{r_\ell}(p)} w_\ell^q\,d\mu,
\quad \tilde{S}_\ell:=R^{q-n}S_\ell,
\end{equation*}
where $q=\frac{2n}{n-2}$.
Then $r_1=R$ and $r_\ell\searrow \frac{R}{2}$, $k_1=0$ and $k_\ell\nearrow \frac{d}{R}$, $\Omega_\ell\subset\Omega_{\ell-1}$ and $w_\ell\leq w_{\ell-1}$ on the set $\Omega_\ell$, and $S_\ell\leq S_{\ell-1}$. Clearly, 
\begin{equation}\label{C-20}
1-\frac{k_\ell}{u}\geq \frac{k_{\ell+1}-k_\ell}{k_{\ell+1}}\geq 2^{-\ell},\quad \mbox{on}\quad \Omega_{\ell+1}.
\end{equation}
Choose $\zeta_\ell\in C^1_c(B_{r_\ell}(p))$ such that
$\zeta_\ell=1$ on $B_{r_{\ell+1}}(p)$ and $|\nabla \zeta_\ell|\leq C \frac{ 2^\ell}{R}$. Assume that  $R\sigma=\lambda\leq d$, so that $\sigma\leq \frac{d}{R}$, $k_\ell+\sigma \leq \frac{2d}{R}$. It follows from \eqref{C-11} and \eqref{C-20} that
\begin{eqnarray}\label{C-21}
\int_{\{u>k_{\ell+1}\}\cap B_{r_\ell}(p)}\!\!\!\!\!\!|\nabla u|^2\zeta_\ell^2\,d\mu \leq n\,2^\ell\!\!\!\int_{\Omega_{\ell}} w_\ell^2|\nabla \zeta_\ell|^2\,d\mu
+2^\ell C_{n,\alpha}\sum\limits^3_{j=0}(k_\ell+\sigma)^{4-j}\!\!\!\int_{\Omega_\ell}w_\ell^j\zeta_\ell^2\,d\mu.
\end{eqnarray}
Moreover,
\begin{eqnarray*}
|\nabla(w_{\ell+1}\zeta_\ell)|^2\leq 2|\nabla u|^2\zeta_\ell^2\,\chi_{\Omega_{\ell+1}}+2w_{\ell+1}^2|\nabla\zeta_\ell|^2.
\end{eqnarray*}
Combining this with \eqref{C-21} and the bound $k_\ell+\sigma\leq\frac{2d}{R}$, we obtain
\begin{eqnarray}\label{C-22}
\int_{\Omega_{\ell}}|\nabla (w_{\ell+1}\zeta_\ell)|^2\,d\mu &\leq& 2\int_{\Omega_\ell}|\nabla u|^2\zeta_\ell^2\,\chi_{\{u>k_{\ell+1}\}}\,d\mu+2\int_{\Omega_\ell}w_{\ell+1}^2|\nabla\zeta_\ell|^2\,d\mu\nonumber\\
&\leq&C_{n,\alpha}B_{n,\alpha}^\ell\Big[\frac{4^\ell}{R^2}\int_{\Omega_\ell}w_\ell^2\,d\mu+
\sum\limits_{j=0}^3(\frac{d}{R})^{4-j}\int_{\Omega_\ell}w_\ell^j\,d\mu\Big],
\end{eqnarray}
where $C_{n,\alpha}$, $B_{n,\alpha}>1$ denote constants depending only on $n$ and $\alpha$, and whose values may vary from line to line. Taking $f=w_{\ell+1}\zeta_\ell$ in \eqref{C-16}, we obtain
\begin{equation}\label{C-23}
S_{\ell+1}^{2/q}\leq C_{n,\alpha}B_{n,\alpha}^\ell\Big[\frac{4^\ell}{R^2}\int_{\Omega_\ell}w_\ell^2\,d\mu+
\sum\limits_{j=0}^3(\frac{d}{R})^{4-j}\int_{\Omega_\ell}w_\ell^j\,d\mu\Big].
\end{equation}

We next estimate the terms on the right-hand side of \eqref{C-23}. Since $w_{\ell-1}\geq k_\ell-k_{\ell-1}\geq
\frac{2^{1-\ell}d}{R}$ on $\Omega_\ell$ for $\ell\geq 2$, we have
\begin{equation}\label{C-24}
\mu(\Omega_\ell)\leq\big(\frac{2^{\ell-1}R}{d}\big)^qS_{\ell-1}.
\end{equation}
By H\"older's inequality and the fact that $S_\ell\leq S_{\ell-1}$, for $0\leq j\leq 3$ and $j\leq q$,
\begin{equation}\label{C-25}
\int_{\Omega_\ell}w^j_\ell\,d\mu\leq S_\ell^{\frac{j}{q}}\mu(\Omega_\ell)^{1-\frac{j}{q}}\leq\big(\frac{2^{\ell-1}R}{d}\big)^{q-j}S_{\ell-1}.
\end{equation}
Substituting \eqref{C-24} and \eqref{C-25} into \eqref{C-23}, we obtain
\begin{equation}\label{C-26}
S_{\ell+1}^{\frac{2}{q}}\leq C_{n,\alpha}B_{n,\alpha}^\ell R^{q-4}(d^{2-q}+d^{4-q})S_{\ell-1}.
\end{equation}
Assume in addition that $0<\lambda\leq d\leq 1$. Then $d^{4-q}\leq d^{2-q}$, and \eqref{C-26} reduces to
\begin{equation*}
S_{\ell+1}\leq C_{n,\alpha}B_{n,\alpha}^\ell\, R^{\frac{(q-4)q}{2}}d^{\frac{(2-q)q}{2}}S^{\frac{q}{2}}_{\ell-1},
\end{equation*}
or, equivalently,
\begin{equation}\label{C-27}
d^{-q}\tilde{S}_{\ell+1}\leq C_{n,\alpha}B_{n,\alpha}^\ell(d^{-q}\tilde{S}_{\ell-1})^{\frac{q}{2}},
\end{equation}
where we have used the fact $n-q=\frac{(n-4)q}{2}$. It is worth mentioning that the required exponent condition $3\leq q$ is equivalent to $n\leq 6$. This explains why​ the estimate holds only for $n\leq 6$.

Set 
\begin{equation*}
\beta:=\frac{q}{2},\quad\quad T_m:=d^{-q}\tilde{S}_{2m-1},\quad\quad \delta_0:=C_{n,\alpha}^{-\frac{1}{\beta-1}}B_{n,\alpha}^{-\frac{2\beta}{(\beta-1)^2}}.
\end{equation*}
It follows from \eqref{C-19} that $\tilde{S}_1\leq {C'}^{\frac{q}{2}}_{n,\alpha}\mathfrak{C}(p,2R)^{\frac{q}{2}}$.
Choose $M_{n,\alpha}$ so large that ${C'}^{\frac{q}{2}}_{n,\alpha}M_{n,\alpha}^{-q}\leq \delta_0$, and set
\begin{equation}\label{C-28}
d:=M_{n,\alpha}(\mathfrak{C}(p,2R)^{\frac{1}{2}}+\lambda).
\end{equation}
Then $d\geq \lambda$, and $d\leq 1$ provided that $\mathfrak{C}(p,2R)+\lambda^2\leq \epsilon_0\leq (8 M_{n,\alpha}^2)^{-1}$. It is easy to check that $T_1\leq \delta_0$, and an induction argument yields
\begin{equation}\label{C-29}
T_m\leq C_{n,\alpha}^{-\frac{1}{\beta-1}}B_{n,\alpha}^{-\frac{2m}{\beta-1}-\frac{2}{(\beta-1)^2}}\rightarrow 0,
\end{equation}
as $m\rightarrow \infty$, where we have used $\beta>1$, and $B_{n,\alpha}>1$. It follows that $S_\ell \rightarrow 0$ as $\ell\rightarrow \infty$, and hence
\begin{equation*}
\int_{\Sigma\cap B_{\frac{R}{2}}(p)}(u-\frac{d}{R})^{+^q}\,d\mu\leq\liminf\limits_{\ell\rightarrow \infty}\int_{\Sigma\cap B_{r_\ell}(p)}w_\ell^q\,d\mu=0.
\end{equation*} 
This implies
\begin{equation}\label{C-30}
u\leq \frac{d}{R}\quad\mbox{on}\quad B_{\frac{R}{2}}(p).
\end{equation}
This completes the proof of Theorem \ref{Theorem1.1}.

\hfill$\Box$

\section{Cone approximation}
The purpose of this section is to study the geometric blow-up argument of EDH-stable hypersurfaces. Since EDH-stationary hypersurfaces are not translation invariant, we must restrict attention to admissible ones within a fixed ball in $\mathbb{R}^{n+1}$.

\subsection{Plane approximation} 
Let $\Sigma$ be a smooth, two-sided hypersurface in $\mathbb{R}^{n+1}$ with unit normal vector field $\nu$, and let $e$ be a given unit vector. Set
\begin{equation*}
h:=\langle x-x_0, e\rangle,\quad
v:=\langle \nu,e\rangle,\quad
g:=(1-v^2)^{\frac{1}{2}}.
\end{equation*}
The following estimate is the weighted version of the basic lemma in the cone approximation argument by Bellettini in Section 3 of \cite{[B-2025]}.

\begin{Lemma}\label{Lemma4.1}
Let $\Sigma\subset\mathbb{R}^{n+1}$ be a hypersurface with 
$\partial\Sigma\cap B_{4R}(x_0)=\emptyset$ that is stable in $B_{4R}(x_0)$ for the $\rho$-weighted functional. Set $\tilde{\rho}=\rho-\rho(x_0)$, and write $d\mu_{\tilde{\rho}}=e^{\tilde{\rho}} d\mathcal{H}^n$. Assume that, for some $\Lambda>0$ and $1\geq \delta>0$,
\begin{equation*}
R^{-n}\mathcal{H}^n(\Sigma\cap B_{4R}(x_0))\leq \Lambda,
\end{equation*}
and
\begin{equation*}
|\tilde{\rho}|_{L^\infty(B_{4R}(x_0))}\leq c_0\delta, \quad
R^m|D^m\rho|_{L^\infty(B_{4R}(x_0))}\leq c_0\delta^m,\quad 1\leq m \leq 3.
\end{equation*}
Then there is a constant $C(n,c_0)>0$ such that
\begin{equation}\label{D-1}
R^{-n}\int_{\Sigma\cap B_{3R}(x_0)} g^2 \,d\mu_{\tilde{\rho}}\leq C(n,c_0) R^{-n-2}\int_{\Sigma\cap B_{\frac{7}{2}R}(x_0)} h^2 \,d\mu_{\tilde{\rho}}+C(n,c_0)\delta^2\Lambda,
\end{equation}
\begin{equation}\label{D-2}
R^{2-n}\int_{\Sigma\cap B_{2R}(x_0)} |A|^2\,d\mu_{\tilde{\rho}}\leq C(n,c_0) R^{-n}\int_{\Sigma\cap B_{3R}(x_0)} g^2 \,d\mu_{\tilde{\rho}}+C(n,c_0)\delta^2\Lambda.
\end{equation}
\end{Lemma}

\emph{Proof}. From the identities $\nabla h=e^{\mathsf{T}}$ and $H=\rho_{\nu}$, we obtain
\begin{eqnarray}\label{D-3}
\Delta_\rho h=\Delta h+\langle \nabla\rho, \nabla h\rangle=\langle H\nu, e\rangle+\langle \nabla\rho, e^{\mathsf{T}}\rangle=\langle D\rho, e\rangle.
\end{eqnarray}
Choose $\eta\in C^1_c(B_{\frac{7}{2}R}(x_0))$ with $0\leq\eta\leq 1$, and $\eta=1$ on $B_{3R}(x_0)$ and $|\nabla\eta|\leq\frac{C}{R}$. Since $|\nabla h|^2=g^2$,
integrating by parts on $\Sigma\cap B_{\frac{7}{2}R}(x_0)$ gives
\begin{eqnarray*}
\int_{\Sigma\cap B_{\frac{7}{2}R}(x_0)}h\eta^2\Delta_\rho h\, d\mu_\rho&=&-\int_{\Sigma\cap B_{\frac{7}{2}R}(x_0)}\langle\nabla(h\eta^2),\nabla h\rangle\, d\mu_{\rho}\\
&=&-\int_{\Sigma\cap B_{\frac{7}{2}R}(x_0)}g^2\eta^2\,d\mu_\rho-2\int_{\Sigma\cap B_{\frac{7}{2}R}(x_0)}h\eta\langle\nabla\eta,\nabla h\rangle\,d\mu_\rho.
\end{eqnarray*}
This, together with \eqref{D-3}, gives
\begin{eqnarray*}
\int_{\Sigma\cap B_{\frac{7}{2}R}(x_0)}g^2\eta^2\,d\mu_\rho&\leq&-\int_{\Sigma\cap B_{\frac{7}{2}R}(x_0)}h\eta^2\langle D\rho,e\rangle \,d\mu_\rho-2\int_{\Sigma\cap B_{\frac{7}{2}R}(x_0)}h\eta\langle\nabla\eta,\nabla h\rangle\,d\mu_\rho\\
&\leq&\int_{\Sigma\cap B_{\frac{7}{2}R}(x_0)}|h|\eta^2|D\rho|\,d\mu_\rho+2\int_{\Sigma\cap B_{\frac{7}{2}R}(x_0)}h^2|\nabla\eta|^2\,d\mu_\rho\\
&{}&+\frac{1}{2}\int_{\Sigma\cap B_{\frac{7}{2}R}(x_0)}g^2\eta^2\,d\mu_\rho,
\end{eqnarray*}
where the last step follows from Young's inequality. Absorbing the last term
into the left-hand side and using $\eta\equiv1$ on $B_{3R}(x_0)$, we get
\begin{align*}
\int_{\Sigma\cap B_{3R}(x_0)} g^2\,d\mu_\rho
&=\int_{\Sigma\cap B_{3R}(x_0)} g^2\eta^2\,d\mu_\rho\\
&\leq 2\int_{\Sigma\cap B_{\frac{7}{2}R}(x_0)} |h|\,\eta^2|D\rho|\,d\mu_\rho
  +4\int_{\Sigma\cap B_{\frac{7}{2}R}(x_0)} h^2|\nabla\eta|^2\,d\mu_\rho\\
&\leq \int_{\Sigma\cap B_{\frac{7}{2}R}(x_0)}
   \Big(\frac{h^2}{R^2}+R^2|D\rho|^2\Big)\,d\mu_\rho
  +4\int_{\Sigma\cap B_{\frac{7}{2}R}(x_0)} h^2|\nabla\eta|^2\,d\mu_\rho\\
&\leq C R^{-2}\int_{\Sigma\cap B_{\frac{7}{2}R}(x_0)} h^2\,d\mu_\rho
  +C\delta^2\int_{\Sigma\cap B_{\frac{7}{2}R}(x_0)} d\mu_\rho ,
\end{align*}
where we have again used Young's inequality together with $\eta\leq1$,
and then $|\nabla\eta|\leq \frac{C}{R}$ and $R|D\rho|_{L^\infty(B_{4R}(x_0))}\leq c_0\delta$.
Finally, multiplying both sides by $R^{-n}e^{-\rho(x_0)}$, and using
$|\tilde\rho|_{L^\infty(B_{4R}(x_0))}\leq c_0\delta$ (so that
$e^{\tilde\rho}\leq e^{c_0\delta}$) together with
$R^{-n}\mathcal{H}^n(\Sigma\cap B_{4R}(x_0))\leq\Lambda$, we obtain \eqref{D-1}.

In a local orthonormal frame field, we have
\begin{equation*}
v_j=-h_{ij}\langle e_i,e\rangle,\quad
v_{jk}=-h_{ijk}\langle e_i,e\rangle -h_{ij}h_{ik}v.
\end{equation*}
Using $H_i=\rho_{i\nu}-h_{ij}\rho_j$ and $\langle \nabla\rho,\nabla v\rangle=-\rho_ih_{ij}\langle e_j,e\rangle$, we
obtain
\begin{eqnarray}\label{D-4}
J_\rho v&=&\Delta_\rho v+|A|^2v-\rho_{\nu\nu} v\nonumber\\
&=&\Delta v+\langle \nabla\rho,\nabla v\rangle +|A|^2v-\rho_{\nu\nu}v\nonumber\\
&=&-D^2\rho(e,\nu)=-\rho_{e\nu},
\end{eqnarray}
where $J_\rho$ denotes the Jacobi operator of the weighted area functional
$\int_\Sigma e^\rho\,d\mathcal{H}^n$.
On the set $\{g>0\}$, set $V=e^{\mathsf{T}}/g$, which is a unit tangent
vector field. Then
\begin{equation}\label{D-5}
|\nabla v|^2=g^2|A(V)|^2,\quad
\nabla g=-\frac{v}{g}\nabla v,\quad
|\nabla v|^2+|\nabla g|^2=|A(V)|^2.
\end{equation}
Using the decomposition $A=\mathring{A}+\frac{H}{n}I$, the elementary
inequality $(a+b)^2\leq(1+\varepsilon)a^2+(1+\varepsilon^{-1})b^2$ with
$\varepsilon=\frac{1}{2(n-1)}$, and the bound
$|\mathring{A}(V)|^2\leq\frac{n-1}{n}|\mathring{A}|^2$ valid for every unit
vector $V$, we get
\begin{equation}\label{D-6}
|A(V)|^2\leq (1-\frac{1}{2n})|A|^2+\frac{2n-1}{n^2}H^2.
\end{equation}
From $g^2=1-v^2$ and \eqref{D-4} we obtain
\begin{equation}\label{D-7}
g\Delta_\rho g=Q v^2+\rho_{e\nu}v-(|\nabla v|^2+|\nabla g|^2),
\end{equation}
where $Q=|A|^2-\rho_{\nu\nu}$. 
Choose a new cutoff $\eta\in C^1_c(B_{3R}(x_0))$ with $0\leq\eta\leq1$,
$\eta=1$ on $B_{2R}(x_0)$ and $|\nabla \eta|\leq\frac{C}{R}$. Integrating by parts on
$\Sigma\cap B_{3R}(x_0)$ gives
\begin{equation}\label{D-8}
\int_{\Sigma\cap B_{3R}(x_0)}g\eta^2\Delta_\rho g\,d\mu_\rho=\int_{\Sigma\cap B_{3R}(x_0)}g^2|\nabla\eta|^2\,d\mu_\rho
-\int_{\Sigma\cap B_{3R}(x_0)}|\nabla(g\eta)|^2\,d\mu_\rho.
\end{equation}
Multiplying \eqref{D-7} by $\eta^2$, substituting \eqref{D-8}, and using
the stability inequality, we arrive at
\begin{equation}\label{D-9}
\int_{\Sigma\cap B_{3R}(x_0)}Q\eta^2\,d\mu_\rho\leq \int_{\Sigma\cap B_{3R}(x_0)}(|\nabla v|^2+|\nabla g|^2-\rho_{e\nu}v)\eta^2\,d\mu_\rho+\int_{\Sigma\cap B_{3R}(x_0)}g^2|\nabla\eta|^2\,d\mu_\rho.
\end{equation}
Combining \eqref{D-9} with \eqref{D-5} and \eqref{D-6}, and using
$|v|\leq1$, $|H|=|\rho_\nu|\leq c_0\delta/R$ and
$|\rho_{\nu\nu}|+|\rho_{e\nu}|\leq c_0\delta^2/R^2$, we obtain
\begin{eqnarray}\label{D-10}
\frac{1}{2n}\int_{\Sigma\cap B_{3R}(x_0)}|A|^2\eta^2\,d\mu_\rho&\leq& \frac{2n-1}{n^2}\int_{\Sigma\cap B_{3R}(x_0)}H^2\eta^2\,d\mu_\rho+\int_{\Sigma\cap B_{3R}(x_0)}\rho_{\nu\nu}\eta^2\,d\mu_\rho\nonumber\\
&{}&-\int_{\Sigma\cap B_{3R}(x_0)}\rho_{e\nu}v\eta^2\,d\mu_\rho+\int_{\Sigma\cap B_{3R}(x_0)}g^2|\nabla\eta|^2\,d\mu_{\rho}\nonumber\\
&\leq& C\int_{\Sigma\cap B_{3R}(x_0)}g^2|\nabla\eta|^2\,d\mu_\rho+C\frac{\delta^2}{R^2}\int_{\Sigma\cap B_{3R}(x_0)}d\mu_{\rho}.
\end{eqnarray}
Similarly, multiplying both sides of \eqref{D-10} by $R^{2-n}e^{-\rho(x_0)}$ and then invoking the bounds $|\tilde{\rho}|_{L^\infty(B_{4R}(x_0))}\leq c_0\delta$ and $R^{-n}\mathcal{H}^n(\Sigma\cap B_{4R}(x_0))\leq \Lambda$, we obtain \eqref{D-2}. It remains to justify the preceding argument on the set $\{g=0\}$, where the
vector field $V=e^{\mathsf{T}}/g$ in \eqref{D-5} is not defined. To this end,
replace $g$ by $g_\varepsilon=\sqrt{g^2+\varepsilon^2}$, which is strictly
positive, and repeat the argument above for each $\varepsilon>0$. Then Fatou's lemma and the dominated convergence theorem show that
\eqref{D-2} holds on $\{g=0\}$ as well.

\hfill$\Box$

\begin{Definition}\label{Definition4.2}
Let $L>4$, $x_0\in\mathbb{R}^{n+1}$, and $R_0>0$. A sequence of smooth, two-sided, properly immersed hypersurfaces $\Sigma_j$ is called \emph{admissible} if
\begin{enumerate}
\item $\partial\Sigma_j\cap B_{LR_0}(x_0)=\emptyset$ and $R_0^{-n}\mu_j(\Sigma_j\cap B_{LR_0}(x_0))\leq \Lambda_L$,
where $\mu_j:=\||\Sigma_j|\|$,
\item $\Sigma_j$ is stable for the $\rho_j$-weighted functional,
\item $|\tilde{\rho}_j|_{L^\infty(B_{LR_0}(x_0))}\leq c_0\delta_j$ where $\tilde{\rho}_j:=\rho_j-\rho_j(x_0)$,
\item $R^m_0|D^m\rho_j|_{L^\infty(B_{LR_0}(x_0))}\leq c_0\delta_j^m$, $1\leq m \leq 3$,
\end{enumerate}
where $\delta_j\searrow 0$ as $j\rightarrow \infty$.
\end{Definition}

We have the following plane approximation result for admissible sequences.
\begin{Proposition}\label{Proposition4.3}
Let $(\Sigma_j,\rho_j)$ be admissible in $B_{8R_0}(x_0)$ and $|\Sigma_j|\rightarrow m|P|$ as integral varifolds, where $P$ is an affine hyperplane through $x_0$ and $m\in\mathbb{N}$. Then
\begin{equation}\label{D-11}
R^{2-n}_0\int_{\Sigma_j\cap B_{4R_0}(x_0)}|A_{\Sigma_j}|^2\,d\mu_{\tilde{\rho}_j}\rightarrow 0,\quad
R_0\sup\limits_{\Sigma_j\cap B_{R_0}(x_0)}|A_{\Sigma_j}|\rightarrow 0.
\end{equation}
Moreover, sheet by sheet, the immersions converge to $mP$ in $C^{2,\beta}_{\mathrm{loc}}(B_{R_0}(x_0))$ for every $0<\beta<1$.
\end{Proposition}

\emph{Proof}. We first prove that the bounded mean curvature and varifold vonvergence imply local Hausdorff convergence. That is, for all $r< 8R_0$, we will show that
\begin{equation}\label{D-12}
\sup\limits_{y\in P\cap B_r(x_0)} \mbox{dist}(y,S_j)\rightarrow 0,\quad
\mbox{and} \quad
\sup\limits_{x\in S_j\cap B_r(x_0)} \mbox{dist}(x,P)\rightarrow 0,
\end{equation}
where $S_j=\mbox{spt}\||\Sigma_j|\|$ is the support of varifold $|\Sigma_j|$. Suppose that the first limit in \eqref{D-12} does not holds. Then there exist $\epsilon >0$ and $y_j\in P\cap B_r(x_0)$ such that $\mbox{dist}(y_j, S_j)\geq \epsilon R_0$ and $y_j\rightarrow y \in P\cap \overline{B_r(x_0)}$. So, for sufficiently large $j$, we have $S_j\cap B_{\frac{\epsilon R_0}{2}}(y)=\emptyset$. Take a nonnegative function $\varphi\in C_c(B_{\frac{\epsilon R_0}{2}}(y))$, and $\varphi\neq 0$ on $P$. Then, $S_j\cap B_{\frac{\epsilon R_0}{2}}(y)=\emptyset$ implies $\int \varphi\,d \||\Sigma_j|\|=0$. However, the varifold convergence $|\Sigma_j|\rightarrow m|P|$ gives $\int \varphi\,d \||\Sigma_j|\|\rightarrow m\int_P \varphi \,d\mathcal{H}^n>0$, which is a contradiction. Now suppose that the second limit in \eqref{D-12} does not hold. Then there exists $\epsilon>0$, and $x_j\in S_j\cap B_r(x_0)$ such that $\mbox{dist}(x_j,P)\geq \epsilon R_0$, and $x_j\rightarrow x$. Take $t=c\epsilon R_0$ with $c>0$ sufficiently small so that $B_{2t}(x_j)\subset\subset B_{8R_0}(x_0)$ and $B_{2t}(x_j)\cap P=\emptyset$. Since $x_j\in S_j$ and $\Sigma_j$ are smooth immersions, we have 
\begin{equation*}
\lim\limits_{s\rightarrow 0}\frac{\mu_j(B_s(x_j))}{\omega_n s^n}\geq 1,
\end{equation*}
where $\mu_j=\||\Sigma_j|\|$. From the admissible condition and the monotonicity formula for bounded mean curvature (see Section 5 in \cite{[A-1972]})
\begin{equation}\label{D-13}
\frac{\mu_j(\Sigma_j\cap B_s(x_j))}{s^n}\leq e^{c_n|H|_{L^\infty}t}\frac{\mu_j(\Sigma_j\cap B_t(x_j))}{t^n},
\quad 0<s<t,
\end{equation}
we have 
\begin{equation*}
\mu_j(B_t(x_j))\geq \omega_n e^{-c_n|H|_{L^\infty}t} t^n.
\end{equation*}
This implies $\mu_j(B_t(x_j))\geq c(n,\epsilon) R_0^n$. Since $\mbox{dist}(x,P)\geq \epsilon R_0$, take $\varphi\in C_c(B_{2t}(x))$ with $\varphi=1$ on $B_t(x)$, and $\mbox{spt}\varphi\cap P=\emptyset$. On a fixed sufficiently small ball, the integral $\int \varphi d\||\Sigma_j|\| \geq c(n,\epsilon) R^n_0>0 $, which contradicts
the varifold convergence $\int \varphi d\||\Sigma_j|\| \rightarrow \int_P\varphi d\mathcal{H}^n=0$.

Taking a unit normal vector $e$ to $P$, we obtain $h_j(x)=\langle x-x_0,e\rangle =\pm\, \mbox{dist}(x,P)$. It follows from the Hausdorff convergence \eqref{D-12} that
\begin{equation}\label{D-14}
R_0^{-1}|h_j|_{L^\infty(\Sigma_j\cap B_{7R_0}(x_0))}\rightarrow 0.
\end{equation}
Taking $R=2R_0$ in Lemma \ref{Lemma4.1} with $\delta=2\delta_j$, we obtain from \eqref{D-14}
\begin{equation}\label{D-15}
R^{-n}_0\int_{\Sigma_j\cap B_{6R_0}(x_0)} g^2_j\,d\mu_{\tilde{\rho}_j}\rightarrow 0, \quad
\mbox{and}\quad 
R^{2-n}_0\int_{\Sigma_j\cap B_{4R_0}(x_0)} |A_{\Sigma_j}|^2\,d\mu_{\tilde{\rho}_j}\rightarrow 0.
\end{equation}
It follows from \eqref{D-15} and the admissible conditions that
\begin{equation*}
(2R_0)^{2-n}\int_{B_{4R_0}(x_0)}(|A_{\Sigma_j}|+\frac{\lambda_j}{2R_0})^2\,d\mu_{\tilde{\rho}_j}\rightarrow 0,
\end{equation*}
where $\lambda_j=2\delta_j$. Applying Corollary \ref{Corollary1.2} with $R=2R_0$, we obtain $2R_0\sup\limits_{\Sigma_j\cap B_{R_0}(x_0)}|A_{\Sigma_j}|\rightarrow 0$. This uniform curvature bound gives a uniform graphical radius. The equation $H_j=(\rho_j)_{\nu_j}$ and the admissible conditions then give interior Schauder estimates. The properness separates the local sheets，and varifold multiplicity counts the sheets.

\hfill$\Box$

\subsection{Classical cone approximation}
Let $\mathcal{C}$ be a stationary cone with vertex $x_0$. Its \emph{spine} is defined by 
\begin{equation*}
S(\mathcal{C}):=x_0+\{v\in\mathbb{R}^{n+1}\,|\,\mathcal{C}+v=\mathcal{C}\},
\end{equation*}
which is the maximum plane contained in $\mathcal{C}$.
A stationary cone $\mathcal{C}$ is called \emph{classical} if it is a sum of at least three distinct halfplanes meeting along an $(n-1)$-dimensional plane, with integer multiplicities.

We need the following $L^q$ and tube estimates near a cone spine to get the classical cone approximation.
\begin{Lemma}\label{Lemma4.4}
Let $(\Sigma_j,\rho_j)$ be admissible in $B_{8R_0}(x_0)$ with $0\leq \delta_j\leq \delta_0$, where $\delta_0=\delta(n,c_0)\leq 1$. Let $S$ be a $d$-dimensional subspace, and set
\begin{equation*}
\mathcal{N}_{tR_0}(x_0+S):=\{x\in B_{4R_0}(x_0)\,|\,\mbox{\emph{dist}}(x,x_0+S)<tR_0\}.
\end{equation*}
Then, for $3\leq n\leq 6$ and $q=\frac{2n}{n-2}$,
\begin{equation}\label{D-16}
R_0^{2-n}\int_{\Sigma_j\cap B_{6R_0}(x_0)}|A_{\Sigma_j}|^2\,d\mu_{\tilde{\rho}_j}\leq C,\quad
R_0^{4-n}\big(\int_{\Sigma_j\cap B_{5R_0}(x_0)}|A_{\Sigma_j}|^q\,d\mu_{\tilde{\rho}_j}\big)^{\frac{2}{q}}\leq C,
\end{equation}
and for $0<t\leq\frac{1}{4}$,
\begin{equation}\label{D-17}
R_0^{-n}\mu_j(\Sigma_j\cap \mathcal{N}_{tR_0}(x_0+S))\leq Ct^{n-d}, \quad
R_0^{2-n}\int_{\Sigma_j\cap\mathcal{N}_{tR_0}(x_0+S)}|A_{\Sigma_j}|^2\,d\mu_{\tilde{\rho}_j}\leq C t^{\frac{2(n-d)}{n}},
\end{equation}
where $\mu_j=\||\Sigma_j|\|$ and $C=C(n,c_0,\Lambda_8,\delta_0)$.
\end{Lemma}

\emph{Proof}. Take $\eta\in C^1_c(B_{7R_0})$ with $0\leq \eta\leq 1$, and $\eta=1$ on $B_{6R_0}(x_0)$, $\nabla\eta\leq \frac{C}{R_0}$. The stability inequality and the admissible conditions give the first estimate in \eqref{D-16},
\begin{eqnarray}\label{D-18}
\int_{\Sigma_j\cap B_{6R_0}}|A_{\Sigma_j}|^2\,d\mu_{\tilde{\rho}_j}&\leq& \int_{\Sigma_j\cap B_{7R_0}}|\nabla\eta|^2\,d\mu_{\tilde{\rho}_j}+\int_{\Sigma_j\cap B_{7R_0}}(\rho_j)_{\nu\nu}\eta^2\,d\mu_{\tilde{\rho}_j}\nonumber\\
&\leq&\frac{C}{R_0^2}\mathcal{H}^n(\Sigma_j\cap B_{8R_0})\leq C R_0^{n-2}.
\end{eqnarray}
Then, taking $\zeta\in C^1_c(B_{6R_0})$ with $0\leq \zeta\leq 1$, and $\zeta=1$ on $B_{5R_0}(x_0)$, $\nabla\zeta\leq \frac{C}{R_0}$, and using \eqref{C-12} and \eqref{D-18}, we obtain
\begin{eqnarray}\label{D-19}
\int_{\Sigma_j\cap B_{6R_0}}|\nabla u_j|^2\zeta^2\,d\mu_{\tilde{\rho}_j}&\leq& C\int_{\Sigma_j\cap B_{6R_0}}
(u_j+\frac{\delta_j}{R_0})^2(|\nabla\zeta|^2+\frac{\delta_j^2}{R_0^2}\zeta^2)\,d\mu_{\tilde{\rho}_j}\nonumber\\
&\leq& CR^{n-4}_0,
\end{eqnarray}
where $u_j=|A_{\Sigma_j}|$. Using the weighted Michael-Simon inequality and \eqref{D-19}, we obtain the $L^q$-estimates
\begin{eqnarray*}
\big(\int_{\Sigma_j\cap B_{5R_0}(x_0)}u_j^q\,d\mu_{\tilde{\rho}_j}\big)^{\frac{2}{q}}&\leq& C\int_{\Sigma_j\cap B_{6R_0}(x_0)}\big(|\nabla(u_j\zeta)|^2+\frac{\delta_j^2}{R_0^2}u_j^2\zeta^2\big)\,d\mu_{\tilde{\rho}_j}\\
&\leq& C\int_{\Sigma_j\cap B_{6R_0}(x_0)}\big(|\nabla u_j|^2\zeta^2+u_j^2|\nabla\zeta|^2+\frac{\delta_j^2}{R_0^2}u_j^2\zeta^2\big)\,d\mu_{\tilde{\rho}_j}\\
&\leq&CR_0^{n-4}.
\end{eqnarray*}

Choose at most $c(n)t^{-d}$ balls $B_{tR_0}(x_a)$ with $a\in I_j$ so as to cover $(x_0+S)\cap B_{5R_0}(x_0)$, and
$\mathcal{N}_{tR_0}(x_0+S)\cap B_{4R_0}(x_0)\subset \bigcup\limits_{a\in I_j} B_{2tR_0}(x_a)$. Let $I'_j:=\{a\in I_j\,|\, \Sigma_j\cap B_{2tR_0}(x_a)\neq \emptyset\}$. Choose $y_a\in \Sigma_j\cap B_{2tR_0}(x_a)$. Then $\Sigma_j\cap B_{2tR_0}(x_a)\subset \Sigma_j\cap B_{4tR_0}(y_a)$ for $a\in I'_j$, and $\Sigma_j\cap\mathcal{N}_{tR_0}(x_0+S)\subset \bigcup\limits_{a\in I'_j}(\Sigma_j\cap B_{4tR_0}(y_a))$. Moreover, $B_{R_0}(y_a)\subset B_{8R_0}(x_0)$ implies $\partial\Sigma_j\cap B_{R_0}(y_a)=\emptyset$. So, the monotonicity formula for bounded mean curvature (see Section 5 in \cite{[A-1972]}) gives
\begin{equation*}
\frac{\mu_j(\Sigma_j\cap B_{4tR_0}(y_a))}{(4tR_0)^n}\leq e^{c_n|H_j|_{L^\infty}R_0}\frac{\mu_j(\Sigma_j\cap B_{R_0}(y_a))}{R_0^n}.
\end{equation*}
This, together with the admissible condition, gives
\begin{equation*}
\mu_j(\Sigma_j\cap B_{4tR_0}(y_a))\leq C t^nR_0^n,
\end{equation*}
where $C=C(n,c_0,\delta_0,\Lambda_8)$, and hence
\begin{equation}\label{D-20}
\mu_j(\Sigma_j\cap \mathcal{N}_{tR_0}(x_0+S))\leq \sum\limits_{a\in I'_j}\mu_j(\Sigma_j\cap B_{4tR_0}(y_a))
\leq Ct^{n-d}R_0^n.
\end{equation}
Then \eqref{D-16} and \eqref{D-20} give
\begin{eqnarray*}
\int_{\Sigma_j\cap \mathcal{N}_{tR_0}(x_0+S)}u^2_j\,d\mu_{\tilde{\rho}_j}&\leq& \big(\int_{\Sigma_j\cap \mathcal{N}_{tR_0}(x_0+S)}u^q_j\,d\mu_{\tilde{\rho}_j}\big)^{\frac{2}{q}}\big(\mu_{\tilde{\rho}_j}(\Sigma_j\cap \mathcal{N}_{tR_0}(x_0+S))\big)^{\frac{2}{n}}\\
&\leq&C\big(\int_{\Sigma_j\cap B_{5R_0}(x_0)}u^q_j\,d\mu_{\tilde{\rho}_j}\big)^{\frac{2}{q}}\big(\mu_j(\Sigma_j\cap \mathcal{N}_{tR_0}(x_0+S))\big)^{\frac{2}{n}}\\
&\leq& Ct^{\frac{2(n-d)}{n}}R_0^{n-2},
\end{eqnarray*}
where $C=C(n,c_0,\delta_0,\Lambda_8)$. Here we have used the fact that the admissibility condition implies $\mu_{\tilde{\rho}_j}(E)\leq C(c_0,\delta_0) \mu_j (E)$ for every Borel set $E$.

\hfill$\Box$

We have the following classical cone approximation for admissible sequences.

\begin{Proposition}\label{Proposition4.5}
Let $(\Sigma_j,\rho_j)$ be admissible in $B_{8R_0}(x_0)$, and $|\Sigma_j|\rightarrow \mathcal{C}$ as integral varifolds, where $\mathcal{C}$ is a classical cone with vertex $x_0$. Then, we have
\begin{equation}\label{D-21}
R_0\sup\limits_{\Sigma_j\cap B_{R_0}(x_0)}|A_{\Sigma_j}|\rightarrow 0.
\end{equation}
Moreover, $\mathcal{C}$ is a finite union of complete hyperplanes, and $\Sigma_j$ converge to $\mathcal{C}$ in $C^{2,\beta}_{\mathrm{loc}}(B_{R_0}(x_0))$ for every $0<\beta<1$.
\end{Proposition}

\emph{Proof}. Let $0<t<\frac{1}{4}$. On the set $B_{4R_0}(x_0)\setminus\mathcal{N}_{tR_0}(S(\mathcal{C}))$, the cone is a finite disjoint union of planar pieces. Covering it by finitely many balls and applying Proposition \ref{Proposition4.3} to each ball gives
\begin{equation*}
R_0^{2-n}\int_{\Sigma_j\cap(B_{4R_0}(x_0)\setminus\mathcal{N}_{tR_0}(S(\mathcal{C})))}|A_{\Sigma_j}|^2\,d\mu_{\tilde{\rho}_j}\rightarrow 0.
\end{equation*}
On the other hand, Lemma \ref{Lemma4.4} with $\mbox{dim}S(\mathcal{C})=n-1$ gives
\begin{equation*}
R_0^{2-n}\int_{\Sigma_j\cap \mathcal{N}_{tR_0}(S(\mathcal{C}))}|A_{\Sigma_j}|^2\,d\mu_{\tilde{\rho}_j}\leq Ct^{\frac{2}{n}}.
\end{equation*}
Thus, we obtain
\begin{equation}\label{D-22}
\limsup\limits_{j\rightarrow \infty}R_0^{2-n}\int_{\Sigma_j\cap B_{4R_0}(x_0)}|A_{\Sigma_j}|^2\,d\mu_{\tilde{\rho}_j}=0.
\end{equation}
Therefore, \eqref{D-22}, Corollary \ref{Corollary1.2} with $R=2R_0$ and $\lambda_j=2\delta_j$, and the admissible condition $\delta_j\rightarrow 0$ give \eqref{D-21}.

For any $x\in S(\mathcal{C})\cap B_{\frac{R_0}{2}}(x_0)$, choose $r$ sufficiently small such that $B_{2r}(x)\subset\subset B_{R_0}(x_0)$. So, we have $\sup\limits_{\Sigma_j\cap B_{2r}(x)}|A_{\Sigma_j}|\rightarrow 0.$
Then the uniform graphical radius exists, and the compactness gives that each limiting sheet is part of an affine plane without boundary. Hence, the sheets extend across the spine, so each half-hyperplane on one side of the spine has an opposite continuation with equal multiplicity. The proof of the $C^{2,\beta}_{\mathrm{loc}}$-convergence  is analogous to that of Proposition \ref{Proposition4.3}.

\hfill$\Box$

\subsection{Finite union of planes approximation}

We finally prove the finite union of planes approximation for admissible sequences.

\begin{Proposition}\label{Proposition4.6}
Let $(\Sigma_j,\rho_j)$ be admissible in $B_{8R_0}(x_0)$, and $|\Sigma_j|\rightarrow \mathcal{C}=\sum\limits_{a=1}^Nm_a |P_a|$ as integral varifolds, where $P_a$ are distinct affine hyperplanes through $x_0$ and $m_a\in\mathbb{N}$. Then we have
\begin{equation}\label{D-23}
R_0\sup\limits_{\Sigma_j\cap B_{R_0}(x_0)}|A_{\Sigma_j}|\rightarrow 0.
\end{equation}
Moreover, the convergence holds in $C^{2,\beta}_{\mathrm{loc}}$ for every $0<\beta<1$.
\end{Proposition}

\emph{Proof}. Clearly, the case $N=1$ is Proposition \ref{Proposition4.3}, and the case $\mbox{dim}S(\mathcal{C})=n-1$ is Proposition \ref{Proposition4.5}. We prove by downward induction on $\mbox{dim}S(\mathcal{C})$. Assume the assertion has been proved for every finite union of hyperplanes whose spine dimension is larger than $d$.

Fix $x\in(\mathcal{C}\cap B_{4R_0}(x_0))\setminus S(\mathcal{C})$, choose $r_x>0$ so small that $B_{8r_x}(x)\subset B_{8R_0}(x_0)$ and $P_a\cap B_{8r_x}(x)=\emptyset$ whenever $x\notin P_a$. By the mean-value theorem and the admissible conditions, we have 
\begin{equation*}
|\rho_j-\rho_j(x)|_{L^\infty(B_{8r_x}(x))}\leq 8c_0\frac{r_x}{R_0}\delta_j,
\end{equation*}
\begin{equation*}
r^m_x|D^m\rho_j|_{L^\infty(B_{8r_x}(x))}\leq c_0(\frac{r_x}{R_0}\delta_j)^m,\quad 1\leq m\leq 3.
\end{equation*}
Thus, $(\Sigma_j,\rho_j)$ is admissible in $B_{8r_x}(x)$, and only the planes containing $x$ occur in the local limit. Their common spine contains both $S(\mathcal{C})$ and the line through $x_0$ and $x$. So, its dimension is at
least $d+1$. The induction hypothesis yields smooth convergence near $x$.  Let $0<t<\frac{1}{4}$. After taking a finite cover, we have
\begin{equation*}
R_0^{2-n}\int_{\Sigma_j\cap(B_{4R_0}(x_0)\setminus\mathcal{N}_{tR_0}(S(\mathcal{C})))}|A_{\Sigma_j}|^2\,d\mu_{\tilde{\rho}_j}\rightarrow 0.
\end{equation*}
For the remaining part, Lemma \ref{Lemma4.4} gives
\begin{equation*}
R_0^{2-n}\int_{\Sigma_j\cap \mathcal{N}_{tR_0}(S(\mathcal{C}))}|A_{\Sigma_j}|^2\,d\mu_{\tilde{\rho}_j}\leq Ct^{\frac{2(n-d)}{n}}.
\end{equation*}
A similar argument as Proposition \ref{Proposition4.5}, we obtain \eqref{D-23}. This implies the $C^{2,\beta}_{\mathrm{loc}}$-convergence; see an explanation in Proposition \ref{Proposition4.3}.

\hfill$\Box$

\emph{Remark}.
In Section 6 we will upgrade the $C^{2,\beta}_{\mathrm{loc}}$-convergence to $C^{\infty}_{\mathrm{loc}}$-convergence, using the fact that the weights there admit uniform $C^\infty$ bounds.

\section{Tangent cone analysis}

Let $\Sigma_j$ be smooth, two-sided, properly immersed hypersurfaces, and write $\mu_j=\||\Sigma_j|\|$ for the mass measure of the associated integral varifold $|\Sigma_j|$. We say $\Sigma_j$ has \emph{locally uniform bounds} if
\begin{equation*}
\sup\limits_j\mu_j(K)<\infty,\quad
\sup\limits_j|H_{\Sigma_j}|_{L^\infty(\Sigma_j\cap K)}<\infty,
\end{equation*}
for every $K\subset\subset\mathbb{R}^{n+1}$.

Let $(\Sigma_j,\rho_j)$ be stable for the weighted functional $\int_\Sigma e^{\rho_j}\,d\mathcal{H}^n$. Let $s_\ell\searrow 0$, and define
\begin{equation*}
\Phi_\ell: \mathbb{R}^{n+1}\rightarrow \mathbb{R}^{n+1},\quad
z\mapsto x+\frac{s_\ell}{R_0}z,
\end{equation*}
where $x\in\mathbb{R}^{n+1}$ and $R_0>0$ are fixed. Define
\begin{equation*}
\Sigma_{j,\ell}:=\Phi_\ell^{-1}(\Sigma_j)=\{z\,|\,  x+\frac{s_\ell}{R_0}z\in\Sigma_j\},\quad
\rho_{j,\ell}:=\rho_j\circ\Phi_\ell-\rho_j(x)=\rho_j( x+\frac{s_\ell}{R_0}z)-\rho_j(x).
\end{equation*}
Then, $\Sigma_{j,\ell}$ is stationary and stable for the $\rho_{j,\ell}-$weighted functional, and
\begin{equation*}
|A_{\Sigma_{j,\ell}}|(z)=\frac{s_\ell}{R_0}|A_{\Sigma_j}|( x+\frac{s_\ell}{R_0}z),\quad
R_0^m|D^m\rho_{j,\ell}|=s^m_\ell|D^m\rho_j|,\quad
1\leq m\leq 3.
\end{equation*}
Here and below the unscaled quantities on the right are evaluated at $ x+\frac{s_\ell}{R_0}z$.

\begin{Definition}\label{Dfinition5.1}
The pairs $(\Sigma_{j,\ell},\rho_{j,\ell})$ are said to be \emph{uniformly diagonally admissible} along $s_\ell$ in $B_{8R_0}(0)$ if
there exist constants $c_0$, $\Lambda_8$, $\epsilon_\ell\searrow 0$ and $I_\ell\nearrow \infty$ such that, for every $j\geq I_\ell$,
\begin{enumerate}
\item $\Sigma_{j,\ell}$ is smooth, two-sided, properly immersed, has no boundary in $B_{8R_0}(0)$, and is stable for the $\rho_{j,\ell}$-weighted functional;

\item $R_0^{-n}\||\Sigma_{j,\ell}|\|(B_{8R_0}(0))\leq\Lambda_8$;

\item $|\rho_{j,\ell}|_{L^\infty(B_{8R_0}(0))}\leq c_0\epsilon_\ell,\quad
R^m_0|D^m\rho_{j,\ell}|_{L^\infty(B_{8R_0}(0))}\leq c_0\epsilon_\ell^m, \quad 1\leq m\leq 3$.
\end{enumerate}
\end{Definition}

\begin{Lemma}\label{Lemma5.2}
Let $2\leq n\leq 6$, $R_0>0$, and let $(\Sigma_j,\rho_j)$ be stable for the $\rho_j$-weighted functional. Suppose $\Sigma_j$ has locally uniform bounds near $x$, and $|\Sigma_j|\rightarrow V$ as integral varifolds. Let $\mathcal{C}$ be one tangent cone to $V$ at $x$, produced by $s_\ell\searrow 0$ in the sense that $(\Phi^{-1}_\ell)_{\#}V\rightarrow \mathcal{C}$ as varifolds in $B_{8R_0}(0)$. Assume that $(\Sigma_{j,\ell},\rho_{j,\ell})$ are diagonally admissible along $s_\ell$ in $B_{8R_0}(0)$, and $\mathcal{C}$ is one of a multiplicity plane, a classical cone, or a finite union of hyperplanes. Then there exists $r>0$ such that
\begin{equation}\label{E-1}
\sup\limits_j\sup\limits_{\Sigma_j\cap B_r(x)}|A_{\Sigma_j}|<\infty.
\end{equation}
Consequently, after passing to a subsequence, the immersions converge in $C^{2,\beta}_{\mathrm{loc}}(B_{\frac{r}{2}}(x))$ to a limit consisting of finitely many immersed sheets, each of which converges as a single-valued graph.
\end{Lemma}

\emph{Proof}. Set $V_\ell:=(\Phi^{-1}_\ell)_{\#}V$, $V_{j,\ell}:=(\Phi^{-1}_\ell)_{\#}|\Sigma_j|=|\Sigma_{j,\ell}|$. Since $\Phi^{-1}_\ell$ is a smooth affine map, $|\Sigma_j|\rightarrow V$ implies
\begin{equation}\label{E-2}
V_{j,\ell}\rightarrow V_\ell,\quad as \quad j\rightarrow \infty.
\end{equation}
On the space of varifolds in $B_{8R_0}(0)$ with a uniform local mass bound, define a metric $d_{8R_0}$ by
\begin{equation*}
d_{8R_0}(W_1,W_2):=\sum\limits_{k=1}^\infty 2^{-k}\frac{|W_1(\varphi_k)-W_2(\varphi_k)|}{1+|W_1(\varphi_k)-W_2(\varphi_k)|},
\end{equation*}
where $\{\varphi_k\}$ is a countable dense family in $C_c(G_n(B_{8R_0}(0)))$. For fixed $\ell$, \eqref{E-2} gives a $J_\ell\in\mathbb{N}$ such that $J_\ell>I_\ell$, and 
\begin{equation}\label{E-3}
d_{8R_0}(V_{j,\ell},V_\ell)<\frac{1}{\ell},\quad for \quad j\geq J_\ell.
\end{equation}

Suppose the conclusion were false. Then for every $r>0$, and $J\in\mathbb{N}$,
\begin{equation}\label{E-4}
\sup\limits_{j>J}\sup\limits_{\Sigma_j\cap B_r(x)}|A_{\Sigma_j}|=\infty.
\end{equation}
Using $r=s_\ell$ in \eqref{E-4}, choose $j(\ell)>J_\ell$ and $z_\ell\in\Sigma_{j(\ell)}\cap B_{s_\ell}(x)$ such that
\begin{equation}\label{E-5}
s_\ell|A_{\Sigma_{j(\ell)}}|(z_\ell)\geq\ell.
\end{equation}
The choice of $j(\ell)$ ensures that $(\Sigma_{j(\ell),\ell},\rho_{j(\ell),\ell})$ is admissible in $B_{8R_0}(0)$. It follows from \eqref{E-3} that $d_{8R_0}(|\Sigma_{j(\ell),\ell}|,V_\ell)<\frac{1}{\ell}$. Clearly, $V_\ell\rightarrow\mathcal{C}$ as varifolds implies $d_{8R_0}(V_\ell,\mathcal{C})\rightarrow 0$. Thus,
\begin{eqnarray*}
d_{8R_0}(|\Sigma_{j(\ell),\ell}|,\mathcal{C})&\leq&d_{8R_0}(|\Sigma_{j(\ell),\ell}|,V_\ell)+d_{8R_0}(V_\ell,\mathcal{C})\\
&<&\frac{1}{\ell}+d_{8R_0}(V_\ell,\mathcal{C})\rightarrow 0.
\end{eqnarray*} 
Therefore, $|\Sigma_{j(\ell),\ell}|\rightarrow \mathcal{C}$ as varifolds in $B_{8R_0}(0)$. Using \eqref{E-5}, we obtain
\begin{equation*}
R_0|A_{\Sigma_{j(\ell),\ell}}|(\tilde{z}_\ell)=s_\ell|A_{\Sigma_{j(\ell)}}|(z_\ell)\geq \ell,
\end{equation*}
where $\tilde{z}_\ell=\frac{R_0}{s_\ell}(z_\ell-x)\in B_{R_0}(0)$, which contradicts Proposition \ref{Proposition4.3}, Proposition \ref{Proposition4.5} and 
Proposition \ref{Proposition4.6}. So, we have found $r_0>0$, $J_0\in \mathbb{N}$, $C_0>0$ such that
\begin{equation*}
\sup\limits_{j>J_0}\sup\limits_{\Sigma_j\cap B_{r_0}(x)}|A_{\Sigma_j}|\leq C_0.
\end{equation*}
Then \eqref{E-1} follows from the smoothness and properness of $\Sigma_j$ with $j\leq J_0$. The proof of the $C^{2,\beta}_{\mathrm{loc}}$-convergence is analogous to that of Proposition \ref{Proposition4.3}.

\hfill$\Box$

Let $V$ be an integral varifold. We say $z\in \mbox{Reg}(V)$, the regular part of $V$, if there is an $r>0$ such that
\begin{equation}\label{E-6}
V\llcorner B_r(z)=\sum\limits_{a=1}^{N_z} m_a|\Gamma_a|, \quad m_a\in\mathbb{N},
\end{equation}
where each $\Gamma_a$ is a smooth embedded hypersurface disk without boundary in $B_r(z)$ and $z\in \Gamma_a$. At a regular point the tangent varifold is unique and is the sum of tangent hyperplanes of local sheets
\begin{equation}\label{E-7}
V_z=\sum\limits_{a=1}^{N_z}m_a|T_z\Gamma_a|.
\end{equation}
A tangent cone of a varifold $V$ at $p$ is any limit
\begin{equation*}
\mathcal{C}_p=\lim\limits_{k\rightarrow\infty}(\eta_{p,r_k})_{\#}V,
\end{equation*}
where $\eta_{p,r_k}(y)=\frac{y-p}{r_k}$ and $r_k\searrow 0$. An iterated tangent cone of $V$ at $x$ is obtained from a finite chain $\mathcal{C}^0,\mathcal{C}^1,\ldots \mathcal{C}^{N}$, where $\mathcal{C}^0$ is a tangent cone of $V$ at $x$, and for $a<N$, $p_a\in\mbox{spt}\|\mathcal{C}^a\|\setminus S(\mathcal{C}^a)$ and $\mathcal{C}^{a+1}$ is a tangent cone of $\mathcal{C}^a$ at $p_a$. We should point out that every finite chain can be realized by a single diagonal sequence of the original hypersurfaces.

\begin{Lemma}\label{Lemma5.3}
Let $\mathcal{C}$ be any iterated tangent cone obtained from the diagonal blow-ups described above. Then the regular part of $\mathcal{C}$ is Euclidean minimal and stable.
\end{Lemma}

\emph{Proof}. We realize the iterated tangent cone $\mathcal{C}$ by a diagonal admissible sequence, without loss of generality, still denoted by $(\Sigma_j,\rho_j)$, so that $|\Sigma_j|\rightarrow \mathcal{C}$ as integral varifolds in $B_{8R_0}(0)$. The weighted stationarity gives, for every compactly supported vector field $X$,
\begin{equation*}
\int_{G_n(U)} (\mbox{div}_TX(y)+\langle D\tilde{\rho}_j, X\rangle)e^{\tilde{\rho}_j}\,d|\Sigma_j|(y,T)=0,
\end{equation*}
where $\mbox{spt}X\subset U\subset B_{8R_0}(0)$ and $\tilde{\rho}_j=\rho_j-\rho_j(0)$. As $j\rightarrow \infty$, $D\tilde{\rho}_j\rightarrow 0$ and $e^{\tilde{\rho}_j}\rightarrow 1$ uniformly on the support of $X$, and hence
\begin{equation}\label{E-8}
\delta_X\mathcal{C}=\int_{G_n(U)} \mbox{div}_TX\,d\mathcal{C}(y,T)=0,
\end{equation}
where the compactness of integral varifolds for the admissible sequence has been used. Hence $\mathcal{C}$ is Euclidean stationary in $B_{8R_0}(0)$.

Fix $z\in \mbox{Reg}(\mathcal{C})$. One tangent cone of $\mathcal{C}$ at $z$ is a finite union of complete tangent hyperplanes with integer multiplicity. By Lemma \ref{Lemma5.2}, there exists $r_z>0$ such that
\begin{equation*}
\sup\limits_{j}\sup\limits_{\Sigma_j\cap B_{r_z}(z)}|A_{\Sigma_j}|<\infty.
\end{equation*}
This curvature bound supplies a uniform graphical radius, and properness resolves the local immersion into finitely parametrized sheets. The varifold convergence identifies the limiting sheets and counts coincident sheets by their integer multiplicities. On a smaller ball, every sheet is a graph of $w_j$ satisfying the prescribed mean curvature equation
\begin{equation}\label{E-9}
\mathcal{M}(w_j)=\langle D\tilde{\rho}_j,\nu_j\rangle,
\end{equation}
where $\mathcal{M}$ is the uniformly elliptic graph mean curvature operator. The convergence is $C^{2,\beta}$ on compact subsets, as in Lemma \ref{Lemma5.2}. Letting $j\rightarrow \infty$ in \eqref{E-9}, we obtain $H_\Gamma=0$, where $\Gamma$ is any one of the local limiting sheets.

To prove stability, set
$U=B_{8R_0}(0)\setminus\operatorname{Sing}(\mathcal{C})$. The preceding curvature
bounds are uniform on every compact subset of $U$. Together with the
local mass bounds and the graphical Schauder estimates, the finite-order
compactness theorem for proper immersions gives, after passing to a
subsequence, a limiting
parameter manifold, denoted by $\Gamma$, which parametrizes $\mathcal{C}$ in $U$,
including its multiplicities. The convergence is $C^2$ on compact subsets
with compatible reparametrizations on an exhaustion of $\Gamma$.
In particular, for every relatively compact open set
$\Omega\subset\subset\Gamma$, there are, for all sufficiently large $j$,
embeddings $\Psi_j:\Omega\to\Sigma_j$ of parameter manifolds such that
the immersed maps converge in $C^2$ on compact subsets of $\Omega$.
Thus the induced metrics converge in $C^1$, and the squared norms of the
second fundamental forms converge uniformly there.

Now take an arbitrary $\varphi\in C_c^1(\Gamma)$ and choose
$\Omega\subset\subset\Gamma$ with
$\operatorname{spt}\varphi\subset\subset\Omega$. Define a single test function
on $\Sigma_j$ by
\[
 \varphi_j(\Psi_j(p))=\varphi(p),\qquad p\in\Omega,
 \qquad
 \varphi_j=0\quad\hbox{outside }\Psi_j(\Omega).
\]
Since $\Psi_j$ is an embedding of parameter manifolds and $\varphi$
vanishes near $\partial\Omega$, this gives
$\varphi_j\in C_c^1(\Sigma_j)$. No restriction to a single graphical
chart is imposed on $\operatorname{spt}\varphi$. Applying the weighted
stability to this entire test function, we obtain
\begin{equation}\label{E-10}
 \int_{\Sigma_j}
 \bigl(|A_{\Sigma_j}|^2-(\tilde\rho_j)_{\nu_j\nu_j}\bigr)
 \varphi_j^2e^{\tilde\rho_j}\,d\mathcal H^n
 \le
 \int_{\Sigma_j}|\nabla\varphi_j|^2e^{\tilde\rho_j}\,d\mathcal H^n.
 \end{equation}
All integrals here are taken on the immersed parameter manifolds. Notice
that $D^2\tilde\rho_j\to0$ and $e^{\tilde\rho_j}\to1$ uniformly on the
relevant compact set. Pulling both integrals back by $\Psi_j$ and using
the $C^2$ convergence, we may pass to the limit in
\eqref{E-10}, giving
\[
 \int_\Gamma |A_\Gamma|^2\varphi^2\,d\mathcal H^n
 \le \int_\Gamma|\nabla\varphi|^2\,d\mathcal H^n.
\]
Here multiplicities are already represented on the limiting parameter
manifold. Since $\varphi\in C_c^1(\Gamma)$ was arbitrary, this proves
stability on the whole regular immersed part in $B_{8R_0}(0)$.

Finally, since $\mathcal{C}$ is a cone, any compactly supported test function on its
regular immersed part can be dilated into $B_{8R_0}(0)$. For
$\varphi_t(x)=\varphi(x/t)$, the stability quadratic form is multiplied
by the positive factor $t^{n-2}$. Thus the same stability inequality holds
on the entire regular immersed part of $\mathcal{C}$. Stationarity extends to the
whole cone in the same way. 
\hfill$\Box$

We now can use the dimension reduction skeleton to obtain the structure of iterated tangent cones. 

\begin{Proposition}\label{Proposition5.4}
For $2\leq n \leq 6$, every iterated tangent cone arising from the admissible sequences above is a finite union of complete hyperplanes with integer multiplicity.
\end{Proposition}

\emph{Proof}. Let $S:=S(\mathcal{C})$, $d:=\dim S$ and $k:=n-d$. Then, we have the orthogonal decomposition $\mathbb{R}^{n+1}=S^\perp\oplus S$, which gives the varifold product $\mathcal{C}=\mathcal{C}'\times |S|$, where $\mathcal{C}'$ is an $k$-dimensional integral cone in $S^\perp=\mathbb{R}^{k+1}$. Hence, we have $\|\mathcal{C}\|=\|\mathcal{C}'\|\otimes\mathcal{H}^d\llcorner S$. The maximality of $S$ gives $S(\mathcal{C}')=\{0\}$.

We use descending induction on $d=\dim S(\mathcal{C})$ within the class of all iterated tangent cones, which are realized by admissible diagonal sequences as before. Assume that every iterated tangent cone whose spine has dimension at least $d+1$ is a finite union of complete hyperplanes, and $\dim S(\mathcal{C})=d$.  

We first prove that $\mbox{Sing}(\mathcal{C})\subset S(\mathcal{C})$. For $d=n$, we immediately obtain $\mathcal{C}=m|S|$, a multiplicity hyperplane. Suppose there were $p\in \mbox{Sing}(\mathcal{C})\setminus S(\mathcal{C})$. Choose any tangent cone $\mathcal{C}_p$ of $\mathcal{C}$ at $p$. Clearly, $S(\mathcal{C})\oplus \mathbb{R}p\subset S(\mathcal{C}_p)$, which gives $\dim S(\mathcal{C}_p)\geq d+1$. The induction hypothesis says that $\mathcal{C}_p$ is a finite union of complete hyperplanes with integer multiplicity. Applying Lemma \ref{Lemma5.2} to a diagonal admissible sequence realizing $\mathcal{C}$, with center $p$ and chosen $\mathcal{C}_p$, gives a uniform curvature bound and sheetwise graphical convergence in a fixed neighborhood of $p$. Since the varifold limit is $\mathcal{C}$, this says that $p\in\mbox{Reg}(\mathcal{C})$, which contradicts the choice of $p$, and hence
\begin{equation}\label{E-11}
\mbox{Sing}(\mathcal{C})\subset S(\mathcal{C}).
\end{equation}
So, the transverse cone $\mathcal{C}'$ is smoothly immersed away from its vertex. Choose $\phi\in C^1_c(\mbox{Reg}(\mathcal{C}'))$ and $\chi\in C^1_c(S)$. By using $\phi(z)\chi(\frac{y}{T})$ as a test function in stablility inequality for the regular part of $\mathcal{C}$ and letting $T\rightarrow \infty$, we see that $\mathcal{C}'$ is stable.

If $k=0$, then $d=n$ and $\mathcal{C}=m|S|$ is already a multiplicity hyperplane. For $k=1$, the transverse stationary integral cone is $\mathcal{C}'=\sum\limits_{a=1}^Nm_a|P_{\nu_a}|$, $P_{\nu_a}=\{t\nu_a\,|\,t\geq 0\}$, $m_a\in\mathbb{N}$, Proposition \ref{Proposition4.5} implies $\mathcal{C}'$ is a finite union of complete lines, and hence $\mathcal{C}$ is again a finite union of complete hyperplanes. For $2\leq k\leq 6$, every sheet in $\mathcal{C}'\setminus\{0\}$ is Euclidean minimal and stable by Lemma \ref{Lemma5.3}. Then, Simons' stable cone theorem (\cite{[S-1968]}) implies that every sheet is a $k$-dimensional hyperplane. Taking the product with $S$ makes each one a complete hyperplane in $\mathbb{R}^{n+1}$.

\hfill$\Box$

\section{Proof of Theorem 1.3}

We prove Theorem \ref{Theorem1.3} by contradiction. Assume that no universal constant in \eqref{A-11} exists. Then
there are data $(\Sigma_j,p_j,R_j)$ satisfying (1), (2), (3) in Theorem \ref{Theorem1.3} and points $z_j\in\Sigma_j\cap B_{R_j}(p_j)$ such that
\begin{equation}\label{F-1}
R_j|A_{\Sigma_j}|(z_j)\rightarrow \infty.
\end{equation}
Fix $R_0>0$, and define
\begin{equation*}
F_j:\mathbb{R}^{n+1}\rightarrow \mathbb{R}^{n+1},\quad y\mapsto z=p_j+\frac{R_j}{R_0}y.
\end{equation*}
Set $\widetilde{\Sigma}_j:=F_j^{-1}(\Sigma_j)$, and $y_j:=\frac{R_0}{R_j}(z_j-p_j)$. Then, $y_j\in B_{R_0}(0)$ and
\begin{equation}\label{F-2}
R_0|A_{\widetilde{\Sigma}_j}|(y)=R_j|A_{\Sigma_j}|(p_j+\frac{R_j}{R_0}y),\quad
R_0|A_{\widetilde{\Sigma}_j}|(y_j)=R_j|A_{\Sigma_j}|(z_j)\rightarrow \infty.
\end{equation}
Set $q_j:=\frac{R_0}{R_j}p_j$ and $\tilde{\rho}_j(y):=\alpha\log\frac{|q_j+y|}{|q_j|}$. Then $|q_j|\geq 8R_0$ and 
$(\widetilde{\Sigma}_j,\tilde{\rho}_j)$ is stable for the $\tilde{\rho}_j$-weighted functional.
For $|y|\leq 4R_0$, we have
\begin{equation}\label{F-3}
\frac{1}{2}\leq\frac{q_j+y}{|q_j|}\leq \frac{3}{2},\quad
|H_{\widetilde{\Sigma}_j}|\leq\frac{|\alpha|}{|q_j|-4R_0}\leq \frac{|\alpha|}{4R_0}.
\end{equation}
On the other hand，the area transformation and assumption (3) give the mass bound
\begin{equation}\label{F-4}
R_0^{-n}\tilde{\mu}_j(\widetilde{\Sigma}_j\cap B_{4R_0}(0))=R^{-n}_j\mu_j(\Sigma_j\cap B_{4R_j}(p_j))\leq 2^{|\alpha|}4^n\Lambda,
\end{equation}
where we have used the fact $2^{-|\alpha|}\leq \frac{|z|^\alpha}{|p_j|^\alpha}\leq 2^{|\alpha|}$ on $B_{R_j}(p_j)$.
Applying Allard's compcatness theorem on $B_{4R_0}(0)$, there exist an integral varifold $V$ and a subsequence (still denoted by $\widetilde{\Sigma}_j$) such that $|\widetilde{\Sigma}_j|\rightarrow V$, and $y_j\rightarrow y_0\in \mbox{spt}(V)\cap\overline{B}_{R_0}(0)$.

Let $s_\ell\searrow 0$, and let $\mathcal{C}=\lim\limits_{\ell\rightarrow \infty}(\eta_{y_0,s_\ell})_{\#}V$ be a tangent cone of $V$ at $y_0$, where $\eta_{y_0,s_\ell}(y)=\frac{y-y_0}{s_\ell}$. Set
\begin{equation}\label{F-5}
\widetilde{\Sigma}_{j,\ell}:=\Phi_\ell^{-1}(\widetilde{\Sigma}_j),\quad
\tilde{\rho}_{j,\ell}:=\tilde{\rho}_j(y_0+\frac{s_\ell}{R_0}z)-\tilde{\rho}_j(y_0)=\alpha\log\frac{|q_j+y_0+\frac{s_\ell}{R_0}z|}{|q_j+y_0|},
\end{equation}
where $\Phi_\ell$ is defined in Section 5. For any fixed $L>4$, for sufficiently large $\ell$, we have $\Phi_\ell(B_{4LR_0})\subset B_{4R_0}(0)$. Hence, for $|z|<4LR_0$, we have
\begin{equation*}
|q_j+y_0+\frac{s_\ell}{R_0}z|\geq |q_j+y_0|-4Ls_\ell\geq 7R_0-4Ls_\ell.
\end{equation*}
It is easy to check that
\begin{equation*}
R_0^m|D^m\tilde{\rho}_{j,\ell}|=s_\ell^m|D^m\tilde{\rho}_j|\leq C(m)|\alpha|\delta^m_{j,\ell},\quad
1\leq m\leq 3.
\end{equation*}
where 
\begin{equation*}
\delta_{j,\ell}=\frac{s_\ell}{|q_j+y_0|-4Ls_\ell}\leq \frac{s_\ell}{7R_0-4Ls_\ell}\rightarrow 0.
\end{equation*}
Choose a diagonal admissible sequence to realize $\mathcal{C}$. Proposition \ref{Proposition5.4} says that $\mathcal{C}$ is a finite union of hyperplanes. Applying Lemma 5.2 with target radius $R_0$, it gives a fixed neighborhood of $y_0$ on which the curvature of $\widetilde{\Sigma}_j$ is uniformly bounded. Since $y_j\rightarrow y_0$, this contradicts \eqref{F-2}. We complete the proof of Theorem \ref{Theorem1.3}.

\vspace{0.5cm}

\noindent \textbf{Acknowledgments}.
This work is supported by the National Natural Science Foundation of China (NSFC) (Grant No. 11871445, No.
12301068), the Stable Support Plan for Youth Teams in Basic Research Fields, Chinese Academy of Sciences (CAS) (Grant No. YSBR-001), and the Fundamental Research Funds for the Central Universities. The authors used ChatGPT to assist in understanding the theory of integral varifolds during the preparation of this manuscript. The
authors take full responsibility for the results of this paper.

\end{document}